\documentclass[a4paper,11pt,reqno]{amsart}
\usepackage[utf8]{inputenc}
\usepackage[english]{babel}
\usepackage{amssymb,amsfonts, amsmath, color}
\usepackage[margin=2.7cm]{geometry}
\usepackage{xcolor}
\usepackage{graphicx, color, enumerate}
\usepackage[active]{srcltx}
\usepackage{pgf}
\usepackage{etex}
\usepackage{verbatim}
\usepackage{pgfkeys}
\usepackage{amsmath}
\usepackage{amsfonts}
\usepackage{amssymb}
\usepackage{amsthm}
\usepackage{float}
\usepackage{xcolor}
\usepackage{mathrsfs}
\usepackage{todonotes}
\usepackage{import}
\usepackage{geometry}
\usepackage[colorlinks,linkcolor={blue},citecolor={blue},urlcolor={black}]{hyperref}
\usepackage{fp}

\newtheorem{theorem}{Theorem}[section]
\newtheorem{proposition}[theorem]{Proposition}

\newtheorem{corollary}[theorem]{Corollary}
\newtheorem{lemma}[theorem]{Lemma}
\newtheorem{remark}[theorem]{Remark}
\newtheorem{definition}[theorem]{Definition}
\newtheorem{assumption}[theorem]{Assumption}

\newcommand{\distpO}[1]{{\sf dist}(#1,\partial\Omega)}

\newcommand{\bcl}{\begin{center}}
\newcommand{\ecl}{\end{center}}
\newcommand{\brl}{\begin{right}}
\newcommand{\erl}{\end{right}}
\newcommand{\ben}{\begin{enumerate}}
\newcommand{\een}{\end{enumerate}}
\newcommand{\overliner}{\begin{array}}
\newcommand{\earr}{\end{array}}
\newcommand{\btab}{\begin{tabular}}
\newcommand{\etab}{\end{tabular}}
\newcommand{\bdoc}{\begin{document}}
\newcommand{\edoc}{\end{document}}
\newcommand{\beqy}{\begin{eqnarray}}
\newcommand{\eeqy}{\end{eqnarray}}

\newcommand{\argmin}{\mathop{\mathrm{argmin}}}

\newcommand{\beqi}{\begin{eqnarray*}}
\newcommand{\eeqi}{\end{eqnarray*}}
\newcommand{\bitem}{\begin{itemize}}

\newcommand{\eitem}{\end{itemize}}
\newcommand{\nln}{\newline}
\newcommand{\newt}{\newtheorem}
\newcommand{\pa}{\partial}
\newcommand{\re}{{I\!\!R}}
\newcommand{\Rn}{\R^N}
\newcommand{\mM}{\mathcal{M}}
\newcommand{\bg}{\boldsymbol{\gamma}}
\newcommand{\vv}{\mathbf{v}}
\newcommand{\ce}{\mathcal{CE}}
\newcommand{\eps}{\varepsilon}

\newcommand{\sign}{{\rm sign}}
\newcommand{\map}{\longrightarrow }
\newcommand{\imp}{\Longrightarrow }
\renewcommand{\div}{\nabla\cdot }
\newcommand{\sen}{{\rm sen\,}}
\newcommand{\tg}{{\rm tg\,}}
\newcommand{\arcsen}{{\rm arcsen\,}}
\newcommand{\arctg}{{\rm arctg\,}}
\newcommand{\supp}{{\textsl supp\ }}
\newcommand{\ity}{\int_{-\iy}^{+\iy}}
\newcommand{\limit}{\lim\limits}
\newcommand{\limi}{\limit_{n\to\infty}}
\newcommand{\sumi}{\sum\limits_{n=1}^{\infty}}
\newcommand{\ulu}{\underline u}
\newcommand{\ulw}{\underline w}
\newcommand{\ulz}{\underline z}
\newcommand{\ulv}{\underline v}
\newcommand{\uls}{\underline s}
\newcommand{\olu}{\overline u}
\newcommand{\olv}{\overline v}
\newcommand{\ols}{\overline s}
\newcommand{\ob}{\overline\b}
\newcommand{\ovar}{\overline\var}
\newcommand{\wv}{\widetilde v}
\newcommand{\wu}{\widetilde u}
\newcommand{\ws}{\widetilde s}
\renewcommand{\a }{\alpha }
\renewcommand{\b }{\beta }
\newcommand{\g }{\gamma}
\newcommand{\G }{\Gamma }
\renewcommand{\d }{\delta }

\newcommand{\D }{\Delta }
\newcommand{\e }{\varepsilon }
\newcommand{\z }{\zeta }
\renewcommand{\l }{\lambda }
\renewcommand{\L }{\Lambda }
\newcommand{\m }{\mu }
\newcommand{\n }{\nabla }
\newcommand{\s }{\sigma }
\newcommand{\Sig }{\Sigma }
\renewcommand{\t }{\tau }
\newcommand{\var }{\varphi }
\renewcommand{\o }{\omega }
\renewcommand{\O }{\Omega }
\newcommand{\R}{{\mathbb{R}}}
\newcommand{\bC}{{\bf C}}
\newcommand{\bZ}{{\bf Z}}
\newcommand{\bN}{{\bf N}}
\newcommand{\bQ}{{\bf Q}}
\newcommand{\bK}{{\bf K}}
\newcommand{\bI}{{\bf I}}
\newcommand{\bv}{{\bf v}}
\newcommand{\bV}{{\bf V}}
\newcommand{\id}{{\sf id}}
\newcommand{\cI}{{\mathcal I}}
\newcommand{\cE}{{\mathcal E}}
\newcommand{\cF}{{\mathcal F}}
\newcommand{\Leb}{{\mathsf{Leb}}}
\newcommand{\bw}{{\boldsymbol{w}}}

\DeclareMathOperator{\suppo}{supp} \DeclareMathOperator{\di}{div}
 \newcommand {\ME}[1]{\todo[inline,size=\footnotesize,color=cyan]{\textbf{ME:} #1}}
\newcommand {\GM}[1]{\todo[inline,size=\footnotesize,color=green]{\textbf{GM:} #1}}

\newenvironment{Proof}{\Rmovelastskip\vskip12pt
plus 1pt \noindent\em\rm}{\hfill {\qed \hskip .2cm}}

\title{Gradient flow for a class of diffusion equations with Dirichlet boundary data}

\author{Matthias Erbar}

\address{\hbox{\parbox{5.7in}{\medskip \noindent{Matthias Erbar, \\Fakult\"at f\"ur Mathematik, \\Universit\"at Bielefeld, \\33501, Bielefeld, Germany \\ [3pt] \emph{E-mail address: }{\tt erbar@math.uni-bielefeld.de}}}}}

\author{Giulia Meglioli}

\address{\hbox{\parbox{5.7in}{\medskip \noindent{Giulia Meglioli, \\Fakult\"at f\"ur Mathematik, \\Universit\"at Bielefeld, \\33501, Bielefeld, Germany \\ [3pt] \emph{E-mail address: }{\tt gmeglioli@math.uni-bielefeld.de}}}}}

\date{}

\begin{document}

\keywords{Gradient flows, Dirichlet boundary conditions, Nonlinear diffusion equation}

%\subjclass[2020]{}

\maketitle

\begin{abstract} In this paper we provide a variational characterisation for a class of non-linear evolution equations with constant non-negative Dirichlet boundary conditions on a bounded domain as gradient flows in the space of non-negative measures. The relevant geometry is given by the  modified Wasserstein distance introduced by Figalli and Gigli that allows for a change of mass by letting the boundary act as a reservoir. We give a dynamic formulation of this distance as an action minimisation problem for curves of non-negative measures satisfying a continuity equation in the spirit of Benamou-Brenier. Then we characterise solutions to non-linear diffusion equations with Dirichlet boundary conditions as metric gradient flows of internal energy functionals in the sense of curves of maximal slope.
\end{abstract}

%\subjclass{35K61, 35K67, 35B99, 35B40, 35B51.}

\bigskip
\smallskip
\section{Introduction}

In this paper, we consider non-linear diffusion equations
\begin{equation}\label{eq:NLD-intro}
\begin{cases}
\partial_t\rho = \Delta L_F(\rho) & \text{in }  \Omega\times(0,+\infty) \\
\rho(0,\cdot) = \rho_0 & \text{in } \Omega \\
\rho = \lambda & \text{on } \partial\Omega\times(0,\infty)\;,
\end{cases}
\end{equation}
with constant Dirichlet boundary condition $\lambda\geq 0$  on an open bounded domain $\Omega\subset \mathbb R^d$. Here $L_F$ is a non-linear function which will be specified later.
Our goal is to give a variational characterisation of solutions in terms of gradient flows in the space of measures equipped with a suitable geometry. Starting from the work of Otto \cite{JKO,O2} a huge number of results characterising various evolutionary PDEs with Neumann boundary conditions as gradient flow in the space of probability measures w.r.t. the Wasserstein distance has been obtained in the literature, see e.g. \cite{ABS, AG13, AGS, Sant, Vil09} for an overview.  However, very little is known to date concerning PDEs with other types of boundary conditions. The first results in this direction have been obtained by Figalli and Gigli \cite{FG10} concerning the linear heat equation with constant Dirichlet boundary conditions that we will briefly describe. As solutions do not conserve mass, the Wasserstein geometry is not appropriate for a gradient flow description. Figalli and Gigli have introduced a variant of the Wasserstein distance allowing for varying mass by letting the boundary $\partial \Omega$ act as a reservoir. Consider the set of non-negative (not necessarily finite measures) on $\Omega$
with finite mean squared distance to the boundary
\[
\mathcal M_2(\Omega):=\left\{\mu\in \mathcal M_+(\Omega): \int_\Omega d(\cdot,\partial\Omega)^2 d\mu<\infty\right\}\;.
\]
For $\mu,\nu$ in $\mathcal M_2(\Omega)$  define 
\begin{equation}\label{eq:Wbp-intro}
Wb_2(\mu,\nu)^2:=\inf_{\bg\in Adm(\mu,\nu)} \int_{\overline\Omega\times\overline\Omega}|x-y|^2d\gamma(x,y)\;,
\end{equation}
where $Adm(\mu,\nu)$ is the set of admissible transport plans consisting of all $\bg\in \mathcal M_+(\overline\Omega\times\overline\Omega)$ such that $\pi^1_\#\bg\vert_{\Omega}=\mu$ and $\pi^2_\#\bg\vert_{\Omega}=\nu$. It is shown in \cite{FG10} 
that $Wb_2$ defines a distance on  $\mathcal M_2(\Omega)$ metrising vague convergence of measures, i.e.~convergence in duality with functions in $C_c(\Omega)$, and sharing many properties with the Wasserstein distance.
Figalli and Gigli consider the time-discrete variational approximation scheme for the gradient flow w.r.t. $Wb_2$ of the Boltzmann entropy given by $\mathcal H(\mu)=\int_\Omega\rho\log\rho\,d\Leb_\Omega$ for $\mu=\rho\Leb_\Omega$. Namely, for a time-step $\tau>0$, the scheme consists in successively solving
\[
\mu^\tau_{n+1} = \underset{\mu}{\rm argmin}\; \mathcal H(\mu) +\frac{1}{2\tau}Wb_2(\mu,\mu_n)^2\;.
\]
It is shown in \cite{FG10} that the curves $(\mu^\tau_t)$ obtained by interpolation (e.g. by setting $\mu^\tau_t=\mu^\tau_n$ on $[n\tau,(n+1)\tau)$) converge in $Wb_2$ as the time step $\tau$ goes to zero to a weak solution $\mu_t=\rho_t\Leb_\Omega$ of the linear diffusion equation $\partial_t\rho =\Delta\rho$ with Dirichlet boundary condition $\rho=1$ on $\partial\Omega$. The boundary condition arises from the fact that $r\mapsto r\log r$ uniquely attains its minimum at $r=1$. By suitably titling the functional $\mathcal H$ other constant boundary conditions $\lambda>0$ can achieved.  These results have been generalised recently by Kim, Koo, and Seo \cite{KKS22} to the porous medium equation $\partial_t\rho=\Delta\rho^\alpha$ by considering the JKO scheme for (suitable tilts of) the functional $\mathcal E_\alpha(\mu)=\int_\Omega \rho^\alpha/(\alpha-1)d\Leb_\Omega$.

This is strong evidence that the heat or porous medium equation should be regarded as the gradient flow of $\mathcal H$ or respectively $\mathcal E_\alpha$ with respect to the distance $Wb_2$. The main contribution of the present paper is to show that indeed the more general class of non-linear diffusion equations with Dirichlet boundary conditions \eqref{eq:NLD-intro} can be characterised as gradient flows of suitable internal energy functionals in the metric space $(\mathcal M_2(\Omega), Wb_2)$ in the framework of gradient flows in metric spaces in the sense of De Giorgi as curves of maximal slope. This question had been left unanswered since \cite{FG10}. The second main contribution of the paper is to give a dynamic characterisation of the transport distance $Wb_2$ in the spirit of the Benamou-Brenier formula \cite{BB2} for the Wasserstein distance. This is also a central ingredient to reach our previous goal. We expect this dynamic point of view to be beneficial also in the analysis of other types of boundary conditions in the future.
\medskip

Let us describe our results in some more detail. Note that the gradient flow of a smooth function $E$ on $\R^n$ or a Riemannian manifold can be characterised as follows: for any smooth curve $(x_t)$ we have  
\begin{equation}\label{eq:gf-smooth}
\frac{d}{dt}E(x_t)=\langle\nabla E(x_t),\dot x_t\rangle\geq  -\frac12|\dot x_t|^2-\frac12|\nabla E(x_t)|^2\;,
\end{equation}
while equality holds if and only if $\dot x_t = -\nabla E(x_t)$, i.e. for the gradient flow curve. The notion of curves of maximal slope generalizes this characterisation to metric spaces, see Section \ref{sec:gradflow} for a recollection of the basic definitions. We consider the metric space $(\mathcal M_2(\Omega), Wb_2)$ and the internal energy functional 
\[
\cF(\mu) = \int_\Omega F(\rho)\,d\Leb_\Omega\;,
\]
for $\mu=\rho\Leb_\Omega$, where $F:[0,\infty)\to[0,\infty)$ is a strictly convex, super-linear function with unique minimum at $\lambda$ such that the non-linearity $L_F$ is given as $$L_F(r):=rF'(r)-F(r),$$ for precise assumptions see Assumption \ref{ass:F}. The role of the modulus of the gradient is taken by the notion of (descending) local slope of $\cF$ given by
\[
|\partial\cF|(\mu):=\limsup_{\nu\to\mu}\frac{\max\{\cF(\mu)-\cF(\nu),0\}}{Wb_2(\mu,\nu)}\;.
\]
The \emph{relaxed slope} $|\partial^-\cF|$ of $\cF$ is defined as the lower semi-continuous envelope of $|\partial\cF|$. The role of the speed of the curve in \eqref{eq:gf-smooth} is taken by the metric derivative of an absolutely continuous curve, see \eqref{eq:def-md}. We show the following result:

\begin{theorem}\label{theo3-intro}
For any absolutely continuous curve $(\mu_t)_{t\in[0,T]}$ in $(\mathcal M_2(\Omega),Wb_2)$ such that $\cF(\mu_0)$ is finite, we have
\[
\mathcal{L}_T(\mu):=\cF(\mu_T)-\cF(\mu_0)+\frac12 \int_0^T\left[|\mu'|^2(r)+|\partial^-\cF|^2(\mu_r)\right]\,dr\ge 0\,.
\]
Moreover, we have $\mathcal{L}_T(\mu_t)=0$ if and only if $\mu_t=\rho_t{\sf Leb}|_\Omega$ such $t\mapsto G(\rho_t)-G(\lambda)$ belongs to $L^2\big([0,T];W^{1,2}_0(\Omega)\big)$ and $(\rho_t)$ is a weak solution to $\partial_t\rho= \Delta f(\rho)$.
\end{theorem}

Here, $G:[0,\infty)\to[0,\infty)$ is the strictly increasing function defined by $$G(0)=0\quad \text{and} \quad G'(r)=\sqrt{r}F''(r).$$ By weak solution we mean a distributional solution, see Definition \ref{def:weak-sol} for details. Note that the condition that $G(\rho_t)$ has trace $G(\lambda)$ on $\partial\Omega$ encodes the Dirichlet boundary condition. In the language of gradient flows in metric spaces, the previous theorem states that $|\partial^-\cF|$ is an upper-gradient of $\cF$ on the metric space $(\mathcal M_2(\Omega),Wb_2)$ and characterises the curves of maximal slope as the solutions to \eqref{eq:NLD-intro}. By applying general results on metric gradient flows \cite{AGS}, we recover the results from \cite{FG10, KKS22} on convergence of the JKO scheme to solutions to \eqref{eq:NLD-intro} as an immediate consequence of Theorem \ref{theo3-intro}.
\medskip

The second main result of this paper is a dynamic characterisation of the transport distance of Figalli-Gigli. We consider slightly more generally the family of distances $Wb_p$ for $p\in[1,\infty)$ defined by \eqref{eq:Wbp-intro} with exponent $p$ instead of $2$. We give a characterisation of absolutely continuous curves w.r.t. the distance $Wb_p$ in terms of solutions to the continuity equation.

\begin{theorem}\label{theo2-intro}
A curve $(\mu_t)_{t\in[0,T]}$ in $(\mathcal M_p(\Omega),Wb_p)$ is absolutely continuous if and only if it is vaguely continuous and there exists a Borel family $(\vv_t)_t$ of vector fields such that
\begin{enumerate}
\item the continuity equation $\partial_t\mu+\nabla\cdot(\mu\vv) =0$ holds in the distributional sense, i.e.
\begin{equation}\label{eq:ce-intro}
\frac{d}{dt}\int\varphi d\mu_t = \int\nabla\varphi \vv_t d\mu_t \qquad \forall \varphi\in C^\infty_c(\Omega)\;,
\end{equation}
\item we have finite $p$-action
\[\int_0^T\|\vv_t\|_{L^p(\mu_t)}dt <\infty\;.\]
\end{enumerate}
In this case, the family of vector fields with minimal $L^p$-norm satisfies $\vert \mu'\vert(t)=||\vv_t||_{L^p(\Omega,\mu_t)}$ for a.e. $t\in [0,T]$, where $|\mu'|$ denotes the metric derivative w.r.t. $Wb_p$.
\end{theorem}

In particular we obtain a dynamic characterisation of the distance $Wb_p$ in the spirit of the Benamou-Brenier formula for the Wasserstein distance \cite{BB2}.
For $\mu_0,\mu_1\in \mathcal M_p(\Omega)$ we have 
\begin{equation*}
Wb_p(\mu_0,\mu_1)^p=\inf\left\{\int_0^1\int_\Omega |\vv_t|^p \,d\mu_t dt\right\}\;,
\end{equation*}
where the infimum is taken over all pairs $(\mu,\vv)$ connecting $\mu_0$ and $\mu_1$ and satisfying \eqref{eq:ce-intro}.

Note that the major difference with the analogous characterisation of absolutely continuous curves w.r.t. the Wasserstein distance (see e.g. \cite[Thm.~8.3.1]{AGS}) is the class of test functions. Requiring \eqref{eq:ce-intro} for all $\varphi\in C^\infty(\Omega)$ prescribes no-flux boundary conditions for the continuity equation and thus conservation of mass. Instead, only requiring \eqref{eq:ce-intro} for $\varphi\in C^\infty_c(\Omega)$ does not fix any boundary conditions and allows for transport to and from the boundary. For the Wasserstein distance any absolutely continuous curve, in particular a gradient flow, necessarily has no-flux boundary conditions. On the other hand, for the distance $Wb_2$, absolutely continous curves can have a variety of boundary conditions. We will see that Dirichlet boundary conditions arise for the gradient flow from the interplay of the distance with the driving functional, more precisely from the finiteness of its slope.
\medskip

We now briefly discuss the proof of Theorem \ref{theo3-intro}. A major challenge stems from the fact that the internal energy $\cF$ typically is not semi-convex along $Wb_2$-geodesics, see \cite[Rmk.~3.4]{FG10}. This is in contrast to the case of the classical Wasserstein distance and prevents the application of general results. Instead, we proceed as follows. We first establish a variational characterisation of solutions to \eqref{eq:NLD-intro} as in Theorem \ref{theo3-intro} but with $|\partial^-\cF|$ replaced by the \emph{energy dissipation functional} $\overline{\mathcal{I}}:\mathcal{M}_2(\Omega)\to[0,+\infty]$, see Theorem \ref{teo2} below. It is defined by first setting
\begin{equation*}
\mathcal{I}(\mu):=
\int_{\Omega}\big|\nabla G(\rho)\big|^2\,d\Leb_\Omega\;,
\end{equation*}
provided that $\mu=\rho\Leb_\Omega$ and $G(\rho)\in W^{1,2}(\Omega)$. Otherwise, we set $\mathcal I(\mu)=+\infty$. Note that for the linear diffusion, i.e. $L_F(r)=r$, we have $G(r)=\sqrt{r}$ and $\mathcal I$ becomes the classical Fisher information. Then, we define $\overline{\mathcal I}(\mu)$ as $\cI(\mu)$ provided $G(\rho)$ additionally has trace $G(\lambda)$ on $\partial \Omega$ and setting $\overline{\mathcal I}(\mu)=+\infty$ else. Since $G$ is strictly increasing this encodes the Dirichlet boundary condition $\lambda$ for $\rho$. This variational description featuring $\overline{\mathcal I}$ is consistent with the fact that the De Giorgi functional $\mathcal L$ of a gradient flow PDE is in many cases strongly related with the path level large deviation rate functional of an underlying particle dynamics, see e.g. \cite{MPR14}. In boundary driven particle systems leading to a macroscopic limit described by a PDE with Dirichlet boundary conditions, the rate function is typically infinite unless the boundary condition is satisfied for all positive times, see e.g. \cite{BDGJL}.

In order to achieve the latter variational description using $\overline{\mathcal I}$, we establish a chain rule for the internal energy (see Proposition \ref{propchainrule} below): for any absolutely continuous curve $(\mu_t)$ in $(\mathcal M_2(\Omega),Wb_2)$ such that $t\mapsto \overline{\mathcal I}(\mu_t)$ is integrable, we have that $t\mapsto \cF(\mu_t)$ is absolutely continuous with
\begin{equation}\label{eq:chain-rule-intro}
\frac{d}{dt}\cF(\mu_t)= \int\langle \bw_t,\vv_t\rangle d\mu_t\;,
\end{equation}
where $(\vv_t)$ is an optimal velocity vector field for $(\mu_t)$ as in Theorem \ref{theo2-intro} and $\bw_t$ is given by $\nabla G(\rho_t)/\sqrt{\rho_t}$. The proof of the chain rule requires a careful regularisation procedure for the curve. From \eqref{eq:chain-rule-intro} we obtain the statement of Theorem \ref{theo3-intro} with $\overline{\mathcal I}$ instead of $|\partial^-\cF|$ immediately via Cauchy-Schwarz and Young inequalities and an analysis of the equality cases.

Finally, in order to obtain the characterisation of solutions to \eqref{eq:NLD-intro} as curves of maximal slope as in Theorem \ref{theo3-intro}, we relate the energy dissipation functional to the slope of $\cF$. Namely, we show in Proposition \ref{prop4} below that 
\[|\partial\cF|^2(\mu)\geq \overline{\mathcal I}(\mu)\qquad \text{for any } \mu\in \mathcal M_2(\Omega)\;.\]
The crucial point here is to show that finiteness of the slope implies that the Dirichlet boundary condition is satisfied, i.e. $\mu=\rho\Leb_\Omega$ with $G(\rho)\in W^{1,2}(\Omega)$ and $G(\rho)$ has trace $G(\lambda)$ on $\partial\Omega$. As $\overline{\mathcal I}$ turns out to be lower semicontinuous one obtains that $\sqrt{\overline{\mathcal I}}$ bounds from below the relaxed slope $|\partial^-\cF|$ as well. Together with Theorem \ref{teo2}, this immediately yields Theorem \ref{theo3-intro}.
\medskip

Finally, we comment on further related results in the literature. During the finalisation of this paper, we learned about independent related work of Quattrocchi \cite{Qu24} and Casteras-Monsaingeon-Santambrogio \cite{CMS}.   In \cite{Qu24}, the author considers linear Fokker-Planck equations on a bounded domain $\Omega$ with general not necessarily constant Dirichlet boundary conditions. He generalises the results in \cite{FG10, Mor18} by proving that solutions can be obtained from a modified scheme of JKO type for a relative entropy where the data are measures supported on the closure $\overline\Omega$ and the role of $Wb_2$ is replaced by a transport-like quantity, which however is not a distance in general. In the special case of constant boundary values in dimension one, he also shows that solutions coincide with curves of maximal slope w.r.t. the relaxed slope. In \cite{CMS} a PDE describing sticky-reflecting diffusions is considered which contains a Dirichlet type compatibility condition for the mass in the domain and on the boundary. The authors show that solutions can be obtained via the JKO scheme, where the relevant functional is sum of bulk and boundary entropies and the distance is the usual Wasserstein distance on the closure of the domain. Moreover, they show that the Dirichlet compatibility condition is encoded via finiteness of the slope.
We also mention the work of Profeta and Sturm \cite{PS20} who give a description of the linear diffusion equation $\partial_t\rho=\Delta\rho$ with homogeneous Dirichet boundary condition $\lambda=0$ as a contraction of a larger auxiliary system of positive and negative densities which can be characterised as a gradient flow. 
\smallskip

\subsection*{Organization of the paper}
In Section 2 we collect preliminaries on the metric space $(\mathcal{M}_p(\Omega),\,Wb_p)$ and on properties of the distance $Wb_p$ for any $p\ge1$. Section 3 is devoted to the continuity equation and the characterisation of absolutely continuous curves in $(\mathcal{M}_p(\Omega),\,Wb_p)$. The proof of this characterisation is given in Section 4. In Section 5 we develop the gradient flow characterisation of non-linear diffusion equations with  Dirichlet boundary conditions.

\subsection*{Acknowledgements}
Funded by the Deutsche Forschungsgemeinschaft (DFG, German Research Foundation) -- Project-ID 317210226 -- SFB 1283. We are grateful for inspiring discussions with Jan Maas and Delio Mugnolo regarding boundary conditions in gradient flows and related topics.

\section{The space $(\mathcal{M}_p(\Omega),\,Wb_p)$}\setcounter{equation}{0}

In this section we recall the definition of the transport distance introduced by Figalli and Gigli in \cite{FG10}. We go slightly beyond the setting considered there by considering modified Wassterstein distances $Wb_p$ with general exponent $p\ge 1$.

Throughout this paper, let  $\Omega\subset \R^d$ be an open bounded domain with Lipschitz boundary. We denote by $\mathcal M_+(\Omega)$ that set of all locally finite Borel measures on $\Omega$ and equip it with the topology of vague convergence, i.e.~convergence in duality with functions in $C_c(\Omega)$. Let us set
\begin{equation}\label{mp}
m_p(\mu):=\int_\Omega d(x,\partial\Omega)^p\,d\mu(x)\;,
\end{equation}
where $d(\cdot,\partial\Omega)$ is the distance from the boundary of $\Omega$.
For $p\in [1,\infty)$ consider the set of measures with finite $p$-th boundary moment, given by
\[
\mathcal{M}_p(\Omega):=\left\{\,\mu\in\mathcal{M}_+(\Omega)\,:\,\,m_p(\mu)<+\infty\right\}\;.
\]
Note that the total mass $\mu(\Omega)$ of $\mu\in \mathcal M_p(\Omega)$ can be infinite.

The following definition, in the case $p=2$, has been given in \cite{FG10}.
\begin{definition}\label{problema}
Let $\mu$, $\nu\in \mathcal{M}_p(\Omega)$. The set of admissible couplings $\textsc{Adm}(\mu,\nu)$ is defined as the set of measures $\boldsymbol{\gamma}$ on $\overline\Omega\times \overline\Omega$ satisfying
\begin{equation}\label{eq01}
\pi_{\sharp}^1 \boldsymbol{\gamma}_{|_{\Omega}}=\mu, \qquad \pi_{\sharp}^2 \boldsymbol{\gamma}_{|_{\Omega}}=\nu.
\end{equation}
For any measure $\boldsymbol{\gamma}$ on $\overline\Omega\times \overline\Omega$, we define its cost $C(\boldsymbol{\gamma})$ as
$$
C(\boldsymbol\gamma):=\int_{\overline{\Omega}\times\overline{\Omega}}|x-y|^p\,\,d \boldsymbol\gamma(x,y)\,.
$$
Then, the distance $Wb_p(\mu,\nu)$ is defined as:
\begin{equation}\label{eqdist}
Wb_p^p(\mu,\nu):=\inf_{\boldsymbol{\gamma}\in \textsc{Adm}(\mu,\nu)} C(\boldsymbol\gamma).
\end{equation}
\end{definition}
The main difference between $Wb_p$ and $W_p$ is the fact that the admissible coupling $\boldsymbol{\gamma}$ is a positive measure on $\overline\Omega\times\overline\Omega$ rather than just $\Omega\times\Omega$ and that the marginals are required to coincide with the given measures only in the interior of $\Omega$. For a more detailed description of such differences see \cite{FG10} in the case of $p=2$.

It is shown in \cite[Prop.~2.7, Prop.~2.9]{FG10},that $(\mathcal{M}_2(\Omega),\,Wb_2)$ is separable, complete and geodesic. The same holds for general $p\geq1$ following the identical arguments.

\subsection{The set of admissible and optimal plans}\label{notation}
Let us observe that, for any $A\subset\Omega$ such that, for some $r>0$, $d(x,\partial\Omega)>r$ for any $x\in A$ and for $\mu\in\mathcal{M}_p(\Omega)$, then
\begin{equation}\label{eq02}
+\infty>m_p(\mu)\ge \int_A r^p\,d\mu(x)=r^p\mu(A),
\end{equation}
i.e. the measure of $A$ is finite.
Let $\boldsymbol{\gamma}\in\mathcal M^+( \overline\Omega\times\overline\Omega)$ be a non-negative measure, we will write $\bg_A^B$ for the restriction of $\bg$ to the rectangle $A\times B\subset \overline\Omega\times\overline\Omega$. Observe that there is a natural splitting of $\bg$ into four parts
$$
\bg=\bg_{\Omega}^{\Omega}+\bg_{\Omega}^{\partial\Omega}+\bg_{\partial\Omega}^{\Omega}+\bg_{\partial\Omega}^{\partial\Omega}.
$$
If $\bg\in\textsc{Adm}(\mu,\nu)$, then
$$
\bg-\bg_{\partial\Omega}^{\partial\Omega}\in\textsc{Adm}(\mu,\nu)\quad \text{and}\quad C(\bg-\bg_{\partial\Omega}^{\partial\Omega})\le C(\bg).
$$
Hence, when looking for optimal plans, it is not restrictive to assume that 
\begin{equation}\label{eq03}
\bg_{\partial\Omega}^{\partial\Omega}=0.
\end{equation}
The set of admissible plans $\textsc{Adm}(\mu,\nu)$ such that \eqref{eq03} is satisfied, is weakly compact, in duality with functions $C_c(\overline\Omega\times\overline\Omega\setminus\partial\Omega\times\partial\Omega)$. Indeed, let $\boldsymbol{\gamma}$ be any measure in $\textsc{Adm}(\mu,\nu)$ and let $A$ be any Borel set in $\Omega$ with positive distance from $\partial\Omega$, we have
$$
\bg(A\times\overline\Omega\cup\overline\Omega\times A)\le \bg(A\times\overline\Omega)+\bg(\overline\Omega\times A)=\mu(A)+\nu(A)<\infty.
$$
Thus from the sequential lower semicontinuity of
$$
\bg \longmapsto C(\bg),
$$
%is weakly lower semicontinuous, i.e. for any $\bg_0\in\textsc{Adm}(\mu,\nu)$ it is true that: for any $\varepsilon>0$ there exists a neighborhood $U(\bg_0)$ such that
%$$
%C(\bg)>C(\bg_0)-\varepsilon\quad \text{for any}\,\,\,\bg\in U(\bg_0).
%$$
%Thus 
we infer the existence of optimal plans. We will denote the set of optimal plans by $\textsc{Opt}(\mu,\nu)$ and we will always assume that an optimal plan satisfies \eqref{eq03}. 

\subsection{Properties of $Wb_p$}
We want to show that $Wb_p$ is a distance over $\mM_p(\Omega)$. To do so, we recall \cite[Lemma 2.1]{FG10}, i.e.
\begin{lemma}\label{lem01}
Let $\mu_1,\mu_2,\mu_3\in \mM_p(\Omega)$ and $\bg^{12}\in\textsc{Adm}(\mu_1,\mu_2)$, $\bg^{23}\in\textsc{Adm}(\mu_2,\mu_3)$ such that $(\bg^{12})_{\partial\Omega}^{\partial\Omega}=(\bg^{23})_{\partial\Omega}^{\partial\Omega}=0$. Then there exists a positive Borel measure $\bg^{123}$ on $\overline\Omega\times\overline\Omega\times\overline\Omega$ such that
\begin{equation}\label{eq04}
\begin{aligned}
&\pi_{\sharp}^{12}\bg^{123}=\bg^{12}+\boldsymbol{\sigma}^{12},\\
&\pi_{\sharp}^{23}\bg^{123}=\bg^{23}+\boldsymbol{\sigma}^{23},
\end{aligned}
\end{equation}
where $\boldsymbol{\sigma}^{12}$ and $\boldsymbol{\sigma}^{23}$ are concentrated on the diagonal of $\partial\Omega\times\partial\Omega$, i.e. on the set of points $\{(x,x):x\in\partial\Omega\}$.
\end{lemma}
\noindent We refer the reader to \cite[Lemma 2.1]{FG10} for the proof of Lemma \ref{lem01}.

\begin{proposition}\label{prop01}
The function $Wb_p$ is a distance on the set $\mM_p(\Omega)$.
\end{proposition}

The proof of Proposition \ref{prop01} is obtained as in \cite[Theorem 2.2]{FG10} with obvious modifications to replace $p=2$ by general $p\ge1$.

\begin{proposition}\label{prop02}
The function $Wb_p$ is lower semicontinuous w.r.t. the weak convergence in duality with functions in $C_c(\Omega)$.
\end{proposition}

\begin{proof}
Let $(\mu_n)_{n\in\mathbb{N}}$ and $(\nu_n)_{n\in\mathbb{N}}$ be two sequence weakly convergent to $\mu$ and $\nu$ respectively. For every $n\in\mathbb{N}$, choose $\bg_n\in\textsc{Opt}(\mu_n,\nu_n)$. Observe that $\bg_n$ is relatively compact in duality with functions in $C_c(\overline\Omega\times\overline\Omega\setminus\partial\Omega\times\partial\Omega)$. Hence, we can estract a subsequence $\bg_{n_k}$ which weakly converge to some $\bg$ in duality with $C_c(\overline\Omega\times\overline\Omega\setminus\partial\Omega\times\partial\Omega)$, i.e.
$$
\int_{\overline\Omega\times\overline\Omega}\varphi\,d\bg_{n_k}\longrightarrow\int_{\overline\Omega\times\overline\Omega}\varphi\,d\bg\quad \text{for any}\,\,\,\varphi\in C_c(\overline\Omega\times\overline\Omega\setminus\partial\Omega\times\partial\Omega).
$$
Now, necessarily 
$$
\pi_{\sharp}^1 \boldsymbol{\gamma}_{|_{\Omega}}=\mu, \qquad \pi_{\sharp}^2 \boldsymbol{\gamma}_{|_{\Omega}}=\nu.
$$
Then,
$$
\begin{aligned}
Wb_p\,^p(\mu,\nu)&\le \int|x-y|^p\,d\bg\le\liminf_{n_k\to+\infty}\int|x-y|^p\,d\bg_{n_k}(x,y)=\liminf_{n_k\to+\infty}Wb_p\,^p(\mu_{n_k},\nu_{n_k}),
\end{aligned}
$$
which ends the proof.

\end{proof}

%\begin{lemma}\label{lem02}
%Let $\mu,\nu\in \mM_p(\Omega)$. For every $1\le p\le q <+\infty$ the following inequality holds
%$$
%Wb_p(\mu,\nu)\le Wb_q(\mu,\nu)\,.
%$$
%\end{lemma}
%\begin{proof}
%We first note that for any $1\le p<q<+\infty$, $\mM_p(\Omega)\subset\mM_q(\Omega)$. Let $\bg\in\textsc{Adm}(\mu,\nu)$, then
%$$
%\left(\int_{\overline\Omega\times\overline\Omega}|x-y|^p\,d\bg\right)^{\frac 1p}=\|x-y\|_{L^p(\bg)}\le\|x-y\|_{L^q(\bg)}=\left(\int_{\overline\Omega\times\overline\Omega}|x-y|^q\,d\bg\right)^{\frac 1q}.
%$$
%Hence $Wb_p\le Wb_q$.
%\end{proof}
%Obviously, from Lemma \ref{lem02}, it follows that $Wb_1\le Wb_p$ for any $p\ge 1$.

%\begin{definition}\label{def2}
%For $p\ge1$, we say that a curve $\mu:(0,T)\to\mM_p(\Omega)$ is $p-$absolutely continuous in $(\mM_p(\Omega),Wb_p)$, if there exists a function $g\in L^p(0,T)$ such that
%\begin{equation}\label{eq15b}
%Wb_p(\mu,\mu_t)\le\int_s^tg(r)\,dr\quad \text{for any}\,\,\,s,t\in(0,T):\,\,s\le t.
%\end{equation}
%\end{definition}
%
%\begin{lemma}\label{lem03}
%Let $p\ge1$. For every $p-$absolutely continuous curve $\mu:(0,T)\to\mM_p(\Omega)$ the metric derivative defined by
%$$
%|\mu'|(t):=\lim_{h\to0}\frac{Wb_p(\mu_t,\mu_{t+h})}{|h|}\,,
%$$
%exists for a.e. $t\in(0,T)$ and $t\to|\mu'|(t)$ belongs to $L^p(0,T)$. The metric derivative $|\mu'|(t)$ is an admissible integrand in the right-hand side of \eqref{eq15b}. Moreover, any other admissible integrand $g\in L^p(0,T)$ satisfies $|\mu'|(t)\le g(t)$ for a.e. $t\in(0,T)$.
%\end{lemma}

\section{Characterisation of absolutly continuous curves}\label{sec3}\setcounter{equation}{0}
In this section, we want to characterise $p-$absolutely continuous curves in the space $(\mM_p(\Omega),Wb_p)$ as solution to the continuity equation
\begin{equation}\label{ce}
\partial_t\mu_t+\nabla\cdot(\mu_t \vv_t) =0,
\end{equation}
with $L^p$-integrable vector fields $\vv_t$.
We start by briefly recalling the notion of absolutely continuous curves in metric spaces. 
\smallskip

A curve $(x_t)_{t\in (a,b)}$ in a complete metric space $(X,d)$ is called $p$-absolutely continuous for
$p\geq1$ if there exists $m\in L^p((a,b))$ such that
\begin{align}\label{eq:abs-continuous}
  d(x_s,x_t)~\leq~\int_s^tm(r) d  r \quad\forall~a\leq s\leq t\leq
  b\ .
\end{align}
In this case we write $x\in AC^p\big((a,b);(X,d)\big)$. For an absolutely continuous curve the metric derivative defined by
\begin{align}\label{eq:def-md}
  |x'_t|~:=~\lim\limits_{h\to0}\frac{d(x_{t+h},x_t)}{|h|}
\end{align}
exists for a.e.~$t$ and is the minimal $m$ in
\eqref{eq:abs-continuous}, see \cite[Thm.1.1.2]{AGS}.
\smallskip

The appropriate notion of weak solution to the continuity equation \eqref{ce} will not specify boundary conditions and will be made precise in the following.
It will be formulated more generally in terms of the vector-valued measures $J_t:=\mu_t\vv_t$ i.e.
\begin{equation}\label{ce2}
\partial_t\mu_t+\nabla\cdot J_t =0.
\end{equation}

\begin{definition}[Continuity equation]\label{def1}
We denote by $\ce_T^\Omega$ the set of all pairs $(\mu, J)$ satisfying the following conditions:
\begin{itemize}
\item[(i)] $\mu:[0,T]\to\mM_p(\Omega)$ is vaguely continuous;
\item[(ii)] $(J_t)_{t\in[0,T]}$ is a Borel family of $\R^d$-valued measures in $\mM(\Omega;\R^d)$;
\item[(iii)] $\int_0^T|J_t|(K) dt < \infty$ for any compact $K\subset\Omega$;
\item[(iv)] for any $\varphi\in C_c^{\infty}(\Omega\times(0,T))$ we have:
\begin{equation}\label{eq:ce-weak}
\int_0^T\left(\int_{\Omega}\partial_t \varphi \,d\mu_t(x)+\int_{\Omega}\nabla\varphi\cdot dJ_t\right)\,dt=0\;.
\end{equation} 
\end{itemize}
Moreover, we will denote by $\ce_T^\Omega(\bar\mu_0,\bar\mu_1)$ the set of pairs $(\mu,J)\in \ce_T^\Omega$ satisfying in addition: $\mu_0=\bar\mu_0$, $\mu_1=\bar\mu_1$. We write $(\mu,\vv)\in \ce_T^\Omega$ if $\vv:[0,T]\times\Omega\to\R^d$ is a time-dependent Borel vector field such that $J_t:=\mu_t\vv_t$  satisfies $(\mu,J)\in \ce_T^\Omega$.
\end{definition}

\begin{remark}\label{rem112}
Observe that the assumption that $\mu_t$ is vaguely continuous is not restrictive. In fact, given a pair $(\mu,J)$ satisfying (ii)-(iv) in the previous definition, there exists a vaguely continuous curve $\tilde\mu:[0,T]\to\mM_p(\Omega)$ such that $\tilde\mu_t=\mu_t$ for a.e. $t\in[0,T]$. See e.g. \cite[Lemma 8.1.2]{AGS} for a proof of the corresponding statement for the continuity equation on $\R^d$.

As a consequence of the vague continuity for any $(\mu,J)\in \ce_T^\Omega$ and $\varphi\in C_c^{\infty}(\Omega\times[0,T])$ we have 
\begin{equation}\label{eq11tris}
\int_0^T\left(\int_{\Omega}\partial_t \varphi \,d\mu_t(x)+\int_{\Omega}\nabla\varphi\cdot dJ_t\right)\,dt-\int_\Omega\varphi(x,T)\,d\mu_T(x)+\int_\Omega\varphi(x,0)\,d\mu_0(x)=0\;.
\end{equation}
\end{remark}

We will need the following representation result for smooth solutions to the continuity equation on a domain without boundary conditions based on the method of characteristics. Consider a time-dependent Borel vector field $\vv:[0,T]\times\overline\Omega\to\R^d$ such that 
\begin{equation}\label{eq:LipVF}
\int_0^T\Big(\sup_{\overline \Omega}|\vv_t|+ {\sf Lip}(\vv_t,\overline{\Omega})\Big) dt <\infty\;.
\end{equation}
Note that for any $x\in\Omega$ and any $s\in[0,T]$ the ODE
\begin{equation}\label{eq:characteristic}
X_s(x,s) =x\;,\quad \frac{d}{dt}X_t(x,s) = \vv_t\big(X_t(x,s)\big)\;,
\end{equation}
admits a unique maximal solution $X_{(\cdot)}(x,s):I(x,s)\to\Omega$ defined in a interval $I(x,s)$ relatively open in $[0,T]$ and containing $s$ as an internal point. Either the solution exists until $t=T$ or at the right endpoint $t_1(x,s)$ of $I(x,s)$ the limit
$\lim_{t\nearrow t_1(x,s)}X_t(x,s)\in \partial\Omega$ exists. Similarly, either the solution extends to $t=0$ or at the left endpoint $t_0(x,s)$ of $I(x,s)$ the limit $\lim_{t\searrow t_0(x,s)}X_t(x,s)\in \partial\Omega$ exists. Note that the solution might reach the boundary $\partial\Omega$ at $t=0$ or $t=T$. Given $s,t\in[0,T]$, we can thus define the map $\Phi_{s,t}:\Omega\to \overline\Omega$ via
\[
\Phi_{s,t}(x):=
\begin{cases}
X_t(x,s)\;, & t\in I(x,s)\;,\\
\lim_{t\nearrow t_1(x,s)}X_t(x,s)\;,& t\geq t_1(x,s)\;,\\
\lim_{t\searrow t_0(x,s)}X_t(x,s)\;,& t\leq  t_0(x,s)\;.
\end{cases}
\]
In other words, the point $\Phi_{s,t}(x)$ is obtained by starting at $x$ at time $s$ and following the characteristics of $\vv$ (forward if $t\geq s$ and backward if $t\leq s$) until time $t$ or until it hits the boundary of $\Omega$. 
Note that $\Phi_{s,t}$ is Lipschitz uniformly in $s,t$.

Let us set $\Omega^{s,t}_{\sf out} := \Phi_{s,t}^{-1}(\partial\Omega)$, the set of points  $x\in\Omega$ such that the characteristic starting at $x$ at time $s$ hits the boundary before or at time $t$. Let $\Omega^{s,t}_{\sf in} := \Phi_{s,t}^{-1}(\Omega)=\Omega\setminus\Omega^{t,s}_{\sf out}$, the set of points $x\in\Omega$ such that the characteristic starting at time $s$ at $x$ exists in $\Omega$ until time $t$. For $t\in(0,T)$ we have 
\[\Omega^{s,t}_{\sf in} =\begin{cases}
\{x\in\Omega : t_1(x,s) > t\} & s\leq t\;,\\
\{x\in\Omega : t_0(x,s) < t\} & s\geq t\;.
\end{cases}
\] 
Note that $\Phi_{s,t}:\Omega^{s,t}_{\sf in}\to\Omega^{t,s}_{\sf in}$ is a bijection , with 
\begin{equation}\label{eq:bijection}
\Phi_{t,s}\circ \Phi_{s,t}=\id \text{ on  }\Omega^{s,t}_{\sf in}\;.
\end{equation}

\begin{proposition}\label{prop:representation-ce}
Let $(\mu,\vv)\in \ce_T^\Omega$ with a Borel vector $\vv:[0,T]\times\overline\Omega\to \R^d$ satisfying \eqref{eq:LipVF}. Then for any $s,t\in [0,T]$ we have
\[ \mu_t\big\vert_{\Omega^{t,s}_{\sf in}} = \big(\Phi_{s,t}\big)_\#\mu_s\big\vert_{\Omega^{s,t}_{\sf in}}\;.\]
Moreover, an admissible coupling $\bg\in\textsc{Adm}(\mu_s,\mu_t)$, according to \eqref{eq03}, is given by,  
\[\begin{aligned}
&\bg_{\Omega}^{\Omega} := (\id,\Phi_{s,t})_\#\mu_s\big\vert_{\Omega^{s,t}_{\sf in}} = (\Phi_{t,s},\id)_\#\mu_t\big\vert_{\Omega^{t,s}_{\sf in}}\;,\\
& \bg_{\Omega}^{\partial\Omega} := (\id,\Phi_{s,t})_\#\mu_s\big\vert_{\Omega^{s,t}_{\sf out}} \;,\\
& \bg_{\partial\Omega}^{\Omega} := (\Phi_{t,s},\id)_\#\mu_t\big\vert_{\Omega^{t,s}_{\sf out}}\;,\\
& \bg_{\partial\Omega}^{\partial\Omega}:=0\;, \end{aligned}\]
and we have the estimate
\begin{equation}\label{eq:cost-est}
Wb_p(\mu_s,\mu_t)^p \leq |t-s|^{p-1}\int_s^t  \int_{\Omega}|\vv_r|^pd\mu_rd r\;.
\end{equation}
\end{proposition}

\begin{proof}
To prove the first claim, 
%one first shows that with $\nu_r:=(\Phi_{s,r})_\# \mu_s\big\vert_{\Omega^{s,t}_{\sf in}}$, the pair $(\nu,\bv)$ is a vaguely continuous solution to the continuity equation \eqref{ce} on $[s,t]$, i.e.~\eqref{eq:ce-weak} holds for any $\varphi\in C^\infty_c\big(\Omega\times(s,t)\big)$. This can be obtained arguing verbatim as in \cite[Lemma 8.1.6]{AGS} where the continuity equation on the whole space is discussed. Then we show that $\nu_r= \mu_t\big\vert_{\Omega^{t,s}_{\sf in}}$. To this end,
let $\psi\in C^\infty_c(\Omega)$ be a test function supported in the open set $\Omega^{t,s}_{\sf in}$ and define $\varphi(r,x):= \psi\big(\Phi_{r,t}(x)\big)$ for $r\in[s,t]$. We deduce from \eqref{eq:bijection}, the continuity of $(r,x)\mapsto \Phi_{t,r}(x)$, and the fact that $\psi$ is supported in $\Omega^{t,s}_{\sf in}$, that $\varphi$ is compactly supported in $[s,t]\times\Omega$ and that $\varphi(s,\cdot)$ is supported in $\Omega^{s,t}_{\sf in}$. Define $\nu_r:=\big(\Phi_{s,r}\big)_\#\mu_s\big\vert_{\Omega^{s,t}_{\sf in}}$. By noting that $\varphi\big(r,\Phi_{s,r}(x)\big)=\psi\big(\Phi_{s,t}(x)\big)$ and differentiating in $r$ we deduce that $\partial_r \varphi+\nabla\varphi\cdot \bv_r=0$. Hence, the continuity equation for $\mu_r$ and the definition of $\nu_r$ entail that 
\[\int \psi d(\mu_t-\nu_t) = \int\varphi(t,\cdot)d(\mu_t-\nu_t) = \int \varphi(s,\cdot)d(\mu_s-\nu_s)=0\;.\]
The arbitrariness of $\psi$ yield that $\mu_t\big\vert_{\Omega^{t,s}_{\sf in}} = \big(\Phi_{s,t}\big)_\#\mu_s\big\vert_{\Omega^{s,t}_{\sf in}}$ as claimed.

It follows immediately that $\bg$ is an admissible coupling and it remains to show \eqref{eq:cost-est}. We estimate the cost of $\bg$ inside the domain as
\begin{align*}
&\int_{\Omega\times\Omega}|x-y|^pd\bg(x,y) 
=
 \int_{\Omega^{s,t}_{\sf in}} |x-\Phi_{s,t}(x)|^pd\mu_s(x)
=
 \int_{\Omega^{s,t}_{\sf in}} \Big|\int_s^t \vv_r\big(\Phi_{s,r}(x)\big)dr\Big|^pd\mu_s(x)\\
 &\leq
 |t-s|^{p-1}\int_s^t\int_{\Omega^{s,t}_{\sf in}}|\vv_r\big(\Phi_{s,r}(x)\big)|^pd\mu_s(x)dr 
 =
 |t-s|^{p-1}\int_s^t\int_{\Phi_{s,r}(\Omega^{s,t}_{\sf in})}|\vv_r|^pd\mu_rdr\;,
\end{align*}
where in the last step we applied the result of the first part of the proof to $[s,r]$ instead of $[s,t]$.
The cost of the part of $\bg$ sending mass to the boundary can be estimated as
\begin{align*}
&\int_{\Omega\times\partial\Omega}|x-y|^pd\bg(x,y) 
=
 \int_{\Omega^{s,t}_{\sf out}} |x-\Phi_{s,t}(x)|^pd\mu_s(x)
=
 \int_{\Omega^{s,t}_{\sf out}} \Big|\int_s^{t_1(x,s)} \vv_r\big(\Phi_{s,r}(x)\big)dr\Big|^pd\mu_s(x)\\
 &\leq
 |t-s|^{p-1}\int_s^t\int_{\Omega^{s,t}_{\sf out}}{\boldsymbol 1}_{\{t_1(x,s)>r\}}|\vv_r\big(\Phi_{s,r}(x)\big)|^pd\mu_s(x)dr \\
& =
 |t-s|^{p-1}\int_s^t\int_{\Omega^{s,t}_{\sf out}\cap \Omega^{s,r}_{\sf in}}|\vv_r\big(\Phi_{s,r}(x)\big)|^pd\mu_s(x)dr =
 |t-s|^{p-1}\int_s^t\int_{\Phi_{s,r}(\Omega^{s,t}_{\sf out}\cap \Omega^{s,r}_{\sf in})}|\vv_r|^pd\mu_rdr\;.
\end{align*}
Similarly, we obtain
\begin{align*}
&\int_{\partial\Omega\times\Omega}|x-y|^pd\bg(x,y) 
\leq 
|t-s|^{p-1}\int_s^t\int_{\Phi_{t,r}(\Omega^{t,s}_{\sf out}\cap \Omega^{t,r}_{\sf in})}|\vv_r|^pd\mu_rdr\;.
\end{align*}
We observe that the three sets
\[
\Phi_{s,r}(\Omega^{s,t}_{\sf in})=\Phi_{t,r}(\Omega^{t,s}_{\sf in})\;,\quad
\Phi_{s,r}(\Omega^{s,t}_{\sf out}\cap \Omega^{s,r}_{\sf in}) = \Omega^{r,t}_{\sf out}\cap\Omega^{r,s}_{\sf in}\;,\quad
\Phi_{t,r}(\Omega^{t,s}_{\sf out}\cap \Omega^{t,r}_{\sf in}) = \Omega^{r,s}_{\sf out}\cap \Omega^{r,t}_{\sf in}
\]
are disjoint. This is immediate upon noting that the first set is contained in $\Omega^{r,t}_{\sf in}\cap\Omega^{r,s}_{\sf in}$. Thus, summing the three estimates above, we obtain \eqref{eq:cost-est}.
\end{proof}

\begin{theorem}\label{teo1}
Let $(\mu_t)_{t\in[0,T]}$ be an absolutely continuous curve in the metric space $\big(\mM_p(\Omega),Wb_p\big)$ and let $|\mu'|\in L^1([0,T])$ be its metric derivative. Then there exists a Borel vector field  $\vv:[0,T]\times\Omega\to\R^d$ such that $(\mu,\vv)\in \ce_T^\Omega$ and
\begin{equation}\label{eq12}
\|\vv_t\|_{L^p(\mu_t,\Omega)}\le |\mu'|(t)\quad \text{ a.e. } t\in [0,T]\,.
\end{equation}
Conversely, if $(\mu,\vv)\in\ce_T^\Omega$ with $\int_0^T\|\vv_t\|_{L^p(\mu_t,\Omega)}\,dt<\infty$, then $(\mu_t)$ is absolutely continuous and
\begin{equation}\label{eq13}
|\mu'|(t)\le \|\vv_t\|_{L^p(\mu_t,\Omega)}\quad \text{a.e.}\,\,\,t\in[0,T].
\end{equation}
\end{theorem}

\begin{proposition}[Benamou-Brenier formula]\label{prop23}
Let $\overline\mu_0,\overline\mu_1\in \mM_p(\Omega)$, then
\begin{equation}\label{eq14}
Wb_p(\overline\mu_0,\overline\mu_1)=\inf\left\{\int_0^T\|\vv_t\|_{L^p(\mu_t,\Omega)}\,dt\;: (\mu,\vv)\in\ce_T^\Omega(\overline\mu_0,\overline\mu_1)\right\}\;.
\end{equation}
\end{proposition}

\begin{remark}\label{rem1}
A standard reparametrization argument %(see e.g.~\cite[]{AGS}) 
shows that we also have
\begin{equation}
Wb_p(\overline\mu_0,\overline\mu_1)^p=\inf\left\{T\int_0^T\|\vv_t\|^p_{L^p(\mu_t,\Omega)}\,dt\,:\,\,(\mu,\vv)\in\ce_T^\Omega(\overline\mu_0,\overline\mu_1)\right\}.
\label{eq15}\end{equation}
\end{remark}

\section{Proof of Theorem \ref{teo1}}\setcounter{equation}{0}
\begin{proof}

\textit{Step $1$:} We show that for any $\varphi\in C_c^{\infty}(\Omega)$ the map $t\mapsto \mu_t(\varphi):=\int_\Omega \varphi d\mu$ is absolutely continuous.\smallskip

Let $s,t\in[0,T]$ and let $\bg\in \textsc{Adm}(\mu_s,\mu_t)$ be an optimal coupling for $Wb_p$. We first observe that, due to Minkowski inequality,
\begin{equation}\label{eqconf}
\begin{aligned}
M_p^{\Omega}(\mu_s)&:=\left(\int_{\overline\Omega}\operatorname{dist}(x,\partial\Omega)^p\,d\mu_s(x)\right)^{1/p}=\left(\int_{\overline\Omega\times \overline\Omega}\operatorname{dist}(x,\partial\Omega)^p\,d\bg(x,y)\right)^{1/p}\\
&\le\left(\int_{\overline\Omega\times \overline\Omega}\left[\operatorname{dist}(y,\partial\Omega)+|x-y|\right]^p\,d\bg(x,y)\right)^{1/p}\\
&\le \left(\int_{\overline\Omega\times \overline\Omega}\operatorname{dist}(y,\partial\Omega)^p\,d\bg(x,y)\right)^{1/p}+\left(\int_{\overline\Omega\times \overline\Omega}|x-y|^p\,d\bg(x,y)\right)^{1/p}\\
&=\left(\int_{\overline\Omega\times \overline\Omega}\operatorname{dist}(y,\partial\Omega)^p\,d\mu_t(y)\right)^{1/p}+Wb_p(\mu_s,\mu_t)\\
&=M_p^\Omega(\mu_t)+Wb_p(\mu_s,\mu_t)\,.
\end{aligned}
\end{equation}
Let $K\subset\subset\Omega$ be the support of $\varphi$ and $\delta:=\operatorname{dist}(K,\partial\Omega)>0$, then we have
$$
\begin{aligned}
|\mu_s(\varphi)-\mu_t(\varphi)|&\le\int_{\overline\Omega\times\overline\Omega}|\varphi(y)-\varphi(x)|\,d\bg\\
%=\left|\int_{\overline\Omega\times\overline\Omega\setminus\Delta}\frac{|\varphi(y)-\varphi(x)|}{|y-x|}\,|y-x|\,d\mu_{st}\right|\\
&\le {\sf Lip}(\varphi)\|\boldsymbol 1_{K\times\bar\Omega\cup \bar\Omega\times K}(x-y)\|_{L^1(\bg,\Omega\times\Omega)}\\
&\le {\sf Lip}(\varphi)\big(\mu_s(K)+\mu_t(K)\big)^{1/q}\|x-y\|_{L^p(\bg,\Omega\times\Omega)}\\
&={\sf Lip}(\varphi)\big(\mu_s(K)+\mu_t(K)\big)^{1/q}\,Wb_p(\mu_s,\mu_t)\,.
\end{aligned}
$$
Observe that, due to \eqref{eqconf}
$$
\begin{aligned}
\mu_t(K)&\le\int_K\frac{\operatorname{dist}(x,\partial\Omega)^p}{\delta^p}\,d\mu_t(x)=\frac{1}{\delta^p}M_p^\Omega(\mu_t)^p\\
&\le\frac{1}{\delta^p}\left(M_p^\Omega(\mu_0)+Wb_p(\mu_0,\mu_t)\right)^p\,.
\end{aligned}
$$
Let $m_p$ be defined as in \eqref{mp}, then we set $C:=\int_0^Tm_p(\mu_s)\,ds$. Finally we get
$$
\begin{aligned}
|\mu_s(\varphi)-\mu_t(\varphi)|&\le\frac{{\sf Lip}(\varphi)}{\delta^{p/q}}\left[\left(M_p^\Omega(\mu_0)+Wb_p(\mu_0,\mu_s)\right)^p+\left(M_p^\Omega(\mu_0)+Wb_p(\mu_0,\mu_t)\right)^p\right]^{1/q}\,Wb_p(\mu_s,\mu_t)\\
&\le\frac{{\sf Lip}(\varphi)}{\delta^{p/q}}\left[\left(M_p^\Omega(\mu_0)+C\right)^p+\left(M_p^\Omega(\mu_0)+C\right)^p\right]^{1/q}\,Wb_p(\mu_s,\mu_t)\,.
\end{aligned}
$$
\smallskip
%\ME{Include estimate on $\mu_s(K)$. Bound in terms of boundary moment, estimate latter by $Wb_p$}

\noindent \textit{Step $2$:} We show that the metric derivative of $\mu_t(\varphi)$ can be estimated with the metric derivative of $\mu_t$.\smallskip

Consider the upper semicontinuous and bounded map
$$
H(x,y):=\begin{cases} \,\,\,\,\,\,\,|\nabla\varphi(x)|\quad &x=y\\ \dfrac{|\varphi(x)-\varphi(y)|}{|x-y|} \quad &x\neq y,\end{cases}
$$
and set $\bg_h\in \textsc{Adm}(\mu_s,\mu_{s+h})$. Then we have
\begin{equation}\label{eq16bis}
\begin{aligned}
\frac{|\mu_{s+h}(\varphi)-\mu_s(\varphi)|}{|h|}
&\le\frac{1}{|h|}\int_{\overline\Omega\times\overline\Omega}|x-y|\,H(x,y)\,d\bg_h\\
&\le\frac{1}{|h|}\left\{\left(\int_{\overline\Omega\times\overline\Omega}|x-y|^p\,d\bg_h\right)^{\frac 1p}\left(\int_{\overline\Omega\times\overline\Omega}H^q(x,y)\,d\bg_h\right)^{\frac1q}\right\}\\
&\le \frac{Wb_p(\mu_{s+h},\mu_s)}{|h|}\left(\int_{\overline\Omega\times\overline\Omega}H^q(x,y)\,d\bg_h\right)^{\frac1q}\,.
\end{aligned}
\end{equation}
Observe that inequality \eqref{eq02} ensures that the family $\{\mu_s\}_s$ is relatively compact in duality with functions in $C_c(\bar\Omega\setminus\partial\Omega)$. This easily implies that also the sequence $\{\bg_h\}$ is relatively compact in duality with functions in $C_c((\bar\Omega\times\bar\Omega)\setminus(\partial\Omega\times\partial\Omega))$. In fact, let $K\subset\left(\bar\Omega\times\bar\Omega\right)\setminus\left(\partial\Omega\times\partial\Omega\right)$ be a compact set and let $0<\delta:={\sf dist}(K,(\partial\Omega\times\partial\Omega))$. Then $K\subset Y:=\{(x,y)\in\bar\Omega\times\bar\Omega\,:\,{\sf dist}((x,y),(\partial\Omega\times\partial\Omega))\ge\delta\}$. Furthermore, we define
$$
\begin{aligned}
& A:=\{x\in\Omega:{\sf dist}(x,\partial\Omega)\ge\delta/\sqrt2\},\\
&B:=\{y\in\Omega:{\sf dist}(y,\partial\Omega)\ge\delta/\sqrt2\}.
\end{aligned}
$$
Let $(x,y)\in K$. We show that $(x,y)\in \left(A\times\bar\Omega\right)\cup\left(\bar\Omega\times B\right)\,.$
By contradiction, if it was not true, then one would have that
$$
{\sf dist}(x,\partial\Omega)<\delta/\sqrt2 \quad \text{and} \quad {\sf dist}(y,\partial\Omega)<\delta/\sqrt2.
$$
The latter implies that $ {\sf dist}((x,y),(\partial\Omega\times\partial\Omega))<\delta$. Therefore we would have $(x,y)\not\in Y$ and $(x,y)\not\in K$ which contradicts the assumption. Hence $K\subset\left(A\times\bar\Omega\right)\cup\left(\bar\Omega\times B\right)$. Then the relative compactness of $\bg_h$ follows from
$$
\bg_h(K)\le\bg_h(A\times\bar\Omega)+\bg_h\left(\bar\Omega\times B\right)=\mu_s(A)+\mu_{s+h}(B)\,.
$$

Therefore, we can extract a subsequence $\{\bg_{h_k}\}$ which converges to some measure $\bg_0$ in duality with functions in $C_c((\bar\Omega\times\bar\Omega)\setminus(\partial\Omega\times\partial\Omega))$. Finally, by \cite[Lemma 2.3]{FG10} we can say that, if $\bg_h\in \textsc{Opt}(\mu_s,\mu_{s+h})$ then $\bg_0$ is optimal as well, i.e. $\bg_0\in \textsc{Opt}(\mu_s,\mu_{s})\equiv(\id,\id)_{\#}\mu_s$. Let $t\in(0,T)$ be a point where $s\mapsto \mu_s$ is metrically differentiable. From \eqref{eq16bis} we get
\begin{equation}\label{eq16}
\begin{aligned}
\limsup_{h\to 0}\frac{|\mu_{s+h}(\varphi)-\mu_s(\varphi)|}{|h|}&\le |\mu'|(t)\left(\int_{\overline\Omega}|H|^q(x,x)\,d\mu_t\right)^{\frac 1q}=|\mu'|(t)\|\nabla\varphi\|_{L^{q}(\Omega,\mu_t)}.
\end{aligned}
\end{equation}
\smallskip

\noindent \textit{Step $3$:} We show that the continuity equation is satisfied in the sense of Definition \ref{def1}.
\smallskip

Set $Q:=\Omega\times(0,T)$ and $\mu:=\int\mu_t\,dt$ be the measure on $Q$ whose disintegration is $\{\mu_t\}_{t}$. Then we have for any $\varphi\in C_c(Q)$ and $h$ small enough
\begin{align*}
\int_Q\partial_s\varphi(x,s)\,d\mu(x,s)&=\lim_{h\to0}\int_Q\frac{\varphi(x,s)-\varphi(x,s-h)}{h}\,d\mu(x,s)\\
&=\lim_{h\to 0}\int_0^T\frac1h\left(\int_{\Omega}\varphi(\cdot,s)\,d\mu_s
-\int_{\Omega}\varphi(\cdot,s)\,d\mu_{s+h}\right)\,ds\;,
\end{align*}
where the last integral is well defined in $(0,T)$ due to the compact support of $\phi$, so that $\phi(\cdot, s)=0$ if $(s+h)\not\in(0,T)$. Hence, using Fatou's lemma, \eqref{eq16}, and H\"older's inequality, we obtain
\begin{align}\nonumber
\left| \int_Q\partial_s\varphi(x,s)\,d\mu(x,s)\right|
&\le\int_0^T\lim_{|h|\to0}\frac{|\mu_s(\varphi)-\mu_{s+h}(\varphi)|}{|h|}\,ds\\\nonumber
&\le\int_0^T|\mu'|(s)\|\nabla\varphi\|_{L^q(\mu_s,\Omega)}\,ds\\ \label{eq17}
&\le\left(\int_0^T|\mu'|^p(s)\,ds\right)^{\frac 1p}\left(\int_{Q}|\nabla\varphi|^q\,d\mu\right)^{\frac 1q}\;.
\end{align}
%where we have used the fact that $d\mu_s(x)ds=d\mu(x,s)$.
Let $V:=\left\{\nabla\varphi:\,\,\varphi\in C_c^{\infty}(\Omega\times(0,T))\right\}$ and denote by $\mathcal V$ its closure in $L^q(\Omega,\mu)$. We define a functional $\mathcal{L}:V\to\R$ via
\[
\mathcal{L}(\nabla\varphi):=-\int_Q\partial_s\varphi(x,s)\,d\mu(x,s)\;,
\]
and observe that $\mathcal{L}$ is linear and bounded by \eqref{eq17} and thus admits a continuous extension to $\mathcal{V}$. Then the problem
\begin{equation}\label{eq18}
\min\left\{\frac 1q\| \mathbf w\|^q_{L^q(\Omega,\mu)}-\mathcal{L}(\mathbf  w)\,\,:\,\, \mathbf w\in \mathcal{V}\right\},
\end{equation}
admits a unique solution $\mathbf{w}\in\mathcal{V}$ and
\[
\vv:=\begin{cases} |\mathbf{w}|^{q-2}\mathbf{w}\quad &\text{if}\,\,\,\mathbf{w}\neq 0\\ 0\quad&\text{if}\,\,\,\mathbf{w}=0,\end{cases}
\]
satisfies
\begin{equation}\label{eq19}
\int_Q\left\langle \vv,\nabla\varphi\right\rangle\,d\mu= \mathcal{L}(\nabla\varphi)\quad \text{for all}\,\,\,\varphi\in C_c^{\infty}\big(\Omega\times(0,T)\big)\;.
\end{equation}
In other words,  $(\mu,\vv)$ satisfies \eqref{eq:ce-weak} and hence $(\mu,\vv)\in \ce_T$. Observe that (iii) in definition \eqref{eq:ce-weak} follows because $\vv_t\in L^p(\mu_t,\Omega)$ and $\mu_t(\Omega)<+\infty$ for a.e. $t\in [0,T]$.
%Set $\vv_t(x):=\vv(x,t)$ and exploit the definition of $\mathcal{L}$ so that we obtain, for every $\varphi\in C_c^{\infty}\big(\Omega\times(0,T)\big)$, that
%\begin{equation}\label{eq110}
%\int_{\Omega}\int_0^1 \left\langle \vv_s(x),\nabla\varphi(x,s)\right\rangle\,d\mu_s(x)\,ds+\int_{\Omega}\int_0^1\partial_s\varphi(x,s)\,d\mu_s(x)\,ds=0\,.
%\end{equation}
\smallskip

\noindent \textit{Step $4$:} We show inequality \eqref{eq12}.
\newline

Let us fix any $J\subset(0,T)$ and $\eta\in C_c^{\infty}(J)$ such that $0\leq \eta\leq 1$. Let $(\nabla\varphi_n)\subset V$ be a sequence converging to $\mathbf{w}$ in $L^q(\mu,Q)$. 
%Then for $\mathit{L}^1$-a.e.~$t\in[0,T]$, there exists a subsequence $n_k$ (possibly depending on $t$) such that $\nabla\varphi_{n_k}\in C_c^{\infty}(\Omega)$ converges in $L^q(\mu_t,\Omega)$ to $\mathbf{w}(\cdot,t)$ where
%$$
%\mathbf{w}(\cdot,t)=\begin{cases} |\vv|^{q-2}\vv\quad &\text{if}\,\,\,\vv\neq 0\\ 0\quad&\text{if}\,\,\,\vv=0.\end{cases}
%$$
Now, due to \eqref{eq17}, \eqref{eq19} we can write
\begin{align*}
\int_Q\eta(s)\,|\vv(x,s)|^p\,d\mu(x,s)&=\int_Q\eta\left\langle \mathbf w,\vv\right\rangle\,d\mu
=\lim_{n\to+\infty}\int_Q\eta \left\langle \vv, \nabla\varphi_n\right\rangle\,d\mu\
=\lim_{n\to+\infty}\mathcal{L}\big(\nabla(\eta\varphi_n)\big)\\
&\le \left(\int_J|\mu'|^p(s)\,ds\right)^{\frac 1p}\,\lim_{n\to+\infty}\left(\int_{\Omega\times J}|\nabla\varphi_n|^q\,d\mu\right)^{\frac 1q}\\
&= \left(\int_J|\mu'|^p(s)\,ds\right)^{\frac 1p}\,\left(\int_{\Omega\times J}|\mathbf{w}|^q\,d\mu\right)^{\frac 1q}\\
&=\left(\int_J|\mu'|^p(s)\,ds\right)^{\frac 1p}\,\left(\int_{\Omega\times J}|\vv|^p\,d\mu\right)^{\frac 1p}.\\
\end{align*}
Letting $\eta$ approximation the characteristic function of $J$, we get
\[
\left(\int_{J}\int_{\Omega}|\vv_s(x)|^p\,d\mu_s(x)\,ds\right)^{1-\frac 1q}\le \left(\int_J|\mu'|^p(s)\,ds\right)^{\frac 1p}\;,
\]
and thus \eqref{eq12}.
\smallskip

\noindent \textit{Step $5$:} We are left to show the converse implication of Theorem \ref{teo1}.

Let $(\mu,J)\in \ce_T^\Omega$ with $J=\vv \mu$ and such that 
\[\int_0^T\|\vv_t\|_{L^p(\mu_t,\Omega)}\,d t <\infty\;.\]

Let $\eta\in C_c^\infty(\R^d)$ be such that $\eta\geq0$, $\operatorname{supp}(\eta)\subset B_1(0)$, and $\|\eta\|_{L^1}=1$. For $\varepsilon>0$ define $\eta_\varepsilon$ by $\eta_\varepsilon(x)=\varepsilon^d\eta(x/\varepsilon)$. Moreover, we set
\[\Omega_\varepsilon = \{x\in \Omega: {\sf dist}(x,\partial\Omega)>\varepsilon\}\;.\]
For each $t$ we define a pair of (vector-valued) measures $(\mu_t^\varepsilon,J_t^\varepsilon)$ on $\Omega_{2\varepsilon}$ by setting
\[\mu_t^\varepsilon := (\mu_t*\eta_\varepsilon)\vert_{\Omega_{2\varepsilon}}\;,
\quad J_t^\varepsilon := (J_t*\eta_\varepsilon)\vert_{\Omega_{2\varepsilon}}\;.\]
Note that the convolutions might not be well-defined on $\R^d$ due to the fact that e.g. $\mu_t$ can have infinite mass in $\Omega$. However, they are well-defined on $\Omega_{2\varepsilon}$. To be more precise, we define $(\mu_t^\varepsilon,J_t^\varepsilon)$ as measures on $\Omega_{2\varepsilon}$ by requiring that for any $\psi\in C_c(\Omega_{2\varepsilon})$, $\Psi\in C_c(\Omega_{2\varepsilon};\R^d)$ we have
\begin{align*}
\int \psi d\mu_t^\varepsilon &= \iint \psi(x) \eta_\varepsilon(x-y)d \mu_t(y)dx\;,\\
\int \Psi dJ_t^\varepsilon &= \iint \eta_\varepsilon(x-y)\Psi(x) d J_t(y)dx\;.
\end{align*}
Note that only $\mu_t\vert_{\Omega_\varepsilon}$ and $J_t\vert_{\Omega_\varepsilon}$ contribute to the above integrals.

A direct computation yields that $(\mu^\varepsilon,J^\varepsilon)\in\ce_T^{\Omega_{2\varepsilon}}$, i.e. the regularised measures satisfy the continuity equation in the smaller domain $\Omega_{2\varepsilon}$. Let us write $J^\varepsilon=\vv^\varepsilon\mu^\varepsilon$. Due to \cite[Lemma 5.20, Proposition 5.21]{Sant} (or also \cite[Lemma 8.1.9, Proposition 8.1.8]{AGS}), we know that $\vv_t^{\varepsilon}$ is Lipschitz and bounded in $x$, uniformly in $t$ and we have the estimate
\begin{equation}
\int_{\Omega_{2\varepsilon}}|\vv_t^{\varepsilon}(x)|^p\,d\mu_t^{\varepsilon}(x)\le\int_{\Omega}|\vv_t|^p\,d\mu_t(x)\quad \text{for every}\,\,\,t\in(0,T)\;.\label{eq113}
\end{equation}
This is based on rewriting
\begin{equation}\label{eq:BBfunct}
\int_{\Omega_{2\varepsilon}}|\vv_t^{\varepsilon}(x)|^p\,d\mu_t^{\varepsilon}(x)=\int_{\Omega_{2\varepsilon}} \alpha_p\Big(\frac{|dJ^\varepsilon}{d\Leb},\frac{d\mu^\varepsilon}{d\Leb}\Big)d\Leb
\end{equation}
with the convex and 1-homogeneous function 
\begin{equation}\label{eq:alpha}
\alpha_p:\R\times [0,\infty)\to [0,\infty]\;,\qquad\alpha_p(v,s) = \begin{cases}
\frac{u^p}{s^{p-1}} & s>0\;,\\
0 & u=s=0\;,\\
+\infty & \text{else .}
\end{cases}
\end{equation}
Then, by Proposition \ref{prop:representation-ce} applied to the domain $\Omega_{2\varepsilon}$ we have
\begin{equation}\label{eq114}
\begin{aligned}
Wb_{p,\Omega_{2\varepsilon}}(\mu^{\varepsilon}_{t},\mu^{\varepsilon}_{(t+h)})
%\left(\int_{\overline\Omega\times\overline\Omega}|x-y|^p\,d\bg^{\varepsilon}\right)^{\frac 1p}\\
%&=\left(\int_{\overline\Omega}\chi_{\{x\in\overline\Omega\,\,s.t.\,\,T_t^{\varepsilon}(x)\in\overline\Omega,\,\,T_{t+h}^{\varepsilon}(x)\in\overline\Omega\}}\,|T_t^{\varepsilon}(x)-T_{t+h}^{\varepsilon}(x)|^p\,d\mu_0^{\varepsilon}\right)^{\frac 1p}\\
%&\le |h|^{\frac 1q}\left(\int_{\overline\Omega}\int_t^{t+h}|\dot{T_s^{\varepsilon}(x)}|^p\,d\mu_s^{\varepsilon}(x)\,ds\right)^{\frac 1p}\\
&\le |h|^{\frac 1q}\left(\int_{\Omega_{2\varepsilon}}\int_t^{t+h}|\vv_s^{\varepsilon}(x)|^p\,d\mu_s^{\varepsilon}(x)\,ds\right)^{\frac 1p}\\
&\le |h|^{\frac 1q}\left(\int_t^{t+h}\int_{\Omega}|\vv_s|^p\,d\mu_s\,ds\right)^{\frac 1p}\\
&=  |h|\left(\frac1{|h|}\int_t^{t+h}\|\vv_s\|^p_{L^p(\mu_s,\Omega)}\,ds\right)^{\frac 1p}.
\end{aligned}
\end{equation}

Here we have made explicit in the notation that the left-hand side is the Figalli-Gigli distance associated to the domain $\Omega_{2\varepsilon}$. To conclude the proof, it is sufficient to  show that for any $s,t\in[0,T]$
\begin{equation}\label{eq:dist-liminf}
Wb_{p,\Omega}(\mu_{s},\mu_{t})\leq \liminf_{\varepsilon\to0}Wb_{p,\Omega_{2\varepsilon}}(\mu^{\varepsilon}_{s},\mu^{\varepsilon}_{t})\;.
\end{equation}
One easily checks that $\mu_t^\varepsilon\to\mu_t$ vaguely for any $t$ as $\varepsilon\to0$.
Hence, the lower semicontinuity of $Wb_{p,\Omega}$ w.r.t. vague convergence \cite[Thm.~2.2]{FG10} yields
\[Wb_{p,\Omega}(\mu_{s},\mu_{t})\leq \liminf_{\varepsilon\to0}Wb_{p,\Omega}(\mu^{\varepsilon}_{s},\mu^{\varepsilon}_{t})\;.\]
It remains to show that 
\begin{equation}\label{eq:dist-comp}
\liminf_{\varepsilon\to0}Wb_{p,\Omega}(\mu^{\varepsilon}_{s},\mu^{\varepsilon}_{t})\leq \liminf_{\varepsilon\to0}Wb_{p,\Omega_{2\varepsilon}}(\mu^{\varepsilon}_{s},\mu^{\varepsilon}_{t})\;.
\end{equation}
To this end, let $\bg^\varepsilon\in\textsc{Adm}_{\Omega_{2\varepsilon}}(\mu^\varepsilon_s,\mu^\varepsilon_t)$ be an optimal admissible transport plan for the domain $\Omega_{2\varepsilon}$ realising $Wb_{p,\Omega_{2\varepsilon}}(\mu^{\varepsilon}_{s},\mu^{\varepsilon}_{t})$.
We construct an admissible plan $\tilde \bg^\varepsilon\in \textsc{Adm}_\Omega(\mu_s^\varepsilon,\mu_t^\varepsilon)$ for the domain $\Omega$ by keeping the transport in the interior of $\Omega_{2\varepsilon}$ and rerouting any mass taken from or sent to a point $x\in \partial\Omega_{2\varepsilon}$ to a point $N(x)\in \partial\Omega$ such that $|x-N(x)|=\distpO{x}$. More precisely, we set
\begin{align*}
\tilde\bg^\varepsilon|_{\Omega_{2\varepsilon}\times\Omega_{2\varepsilon}}&=\bg^\varepsilon|_{\Omega_{2\varepsilon}\times\Omega_{2\varepsilon}}\;,
\quad \tilde\bg^\varepsilon|_{\overline\Omega \times(\overline\Omega\setminus\Omega_{2\varepsilon})}=
 \tilde\bg^\varepsilon|_{(\overline\Omega\setminus\Omega_{2\varepsilon})\times\overline\Omega}=0\;,\\
 \tilde\bg^\varepsilon|_{\partial\Omega\times\Omega_{2\varepsilon}}&=(N,{\id})_\#\bg^\varepsilon|_{\partial\Omega_{2\varepsilon}\times\Omega_{2\varepsilon}}\;,\quad 
\tilde\bg^\varepsilon|_{\Omega_{2\varepsilon}\times\partial\Omega}=({\id,N})_\#\bg^\varepsilon|_{\Omega_{2\varepsilon}\times \partial\Omega_{2\varepsilon}}\;.
\end{align*}
This yields
\begin{align*}
Wb_{p,\Omega}(\mu^{\varepsilon}_{s},\mu^{\varepsilon}_{t})^p
&\leq
\int |x-y|^pd\tilde\bg^\varepsilon\\
&=
\int\limits_{\Omega_{2\varepsilon}\times\Omega_{2\varepsilon}} |x-y|^pd \bg^\varepsilon
+\int\limits_{\partial\Omega_{2\varepsilon}\times\Omega_{2\varepsilon}} |N(x)-y|^pd \bg^\varepsilon + \int\limits_{\Omega_{2\varepsilon}\times\partial\Omega_{2\varepsilon}} |x-N(y)|^pd \bg^\varepsilon\;.
\end{align*}
Note that for any $\delta>0$ there is $C(\delta)>0$ such that $(a+b)^p\leq (1+\delta)a^p + C(\delta) b^p$ for all $a,b\geq0$. 
Combing with the estimate $|N(x)-y|\leq |x-y|+2\varepsilon$, we obtain
\begin{align*}
Wb_{p,\Omega}(\mu^{\varepsilon}_{s},\mu^{\varepsilon}_{t})^p
&\leq
(1+\delta) \cdot Wb_{p,\Omega_{2\varepsilon}}(\mu^{\varepsilon}_{s},\mu^{\varepsilon}_{t})^p
+C(\delta) (2\varepsilon)^p \big(\mu^\varepsilon_s(\Omega_{2\varepsilon})+\mu^\varepsilon_t(\Omega_{2\varepsilon})\big)\;.
\end{align*}
It is sufficent to show that the last term goes to $0$ for $\delta$ fixed as $\varepsilon\to0$. By construction we readily verify that $\mu^\varepsilon_s(\Omega_{2\varepsilon}) \leq \mu_s(\Omega_\varepsilon)$. Since $\mu_s\in \mathcal M_p(\Omega)$, i.e.~$\operatorname{dist}(\cdot,\partial\Omega)^p\in L^1(\mu_s)$, the de la Vall\'ee-Poussin theorem provide non-negative, increasing, and super-linear function $G$ such that $G\big(\operatorname{dist}(\cdot,\partial\Omega)^p\big)\in L^1(\mu_s)$. Then we have by Markov inequality and the super-linearity of $G$ that
\begin{align*}
\varepsilon^p\mu_s(\Omega_\varepsilon) \leq \frac{\varepsilon^p}{G(\varepsilon^p)}\int_{\Omega_\varepsilon} G\big(\operatorname{dist}(\cdot,\partial\Omega)^p\big)d\mu_s(x)\to 0\quad \text{as }\varepsilon\to 0\;.
\end{align*}
The term with $\mu_t^\varepsilon$ can be treated the same way, which finishes the proof.

\end{proof}

\subsection{Proof of Proposition \ref{prop23}}\label{p2}%\setcounter{equation}{0}

\begin{proof}
Let us show the inequality "$\le$". Let $(\mu,\vv)\in\ce^\Omega_T$ and assume
$
\int_0^T\|\vv_t\|_{L^p(\mu_t,\Omega)}\,dt<\infty,
$
as otherwise there is nothing to show. Then, we can apply Theorem \ref{teo1}, therefore, by \eqref{eq13}, we get
$$
Wb_p(\mu_0,\mu_1)\le\int_0^T|\mu'|(t)\,dt\le \int_0^T\|\vv_t\|_{L^p(\mu_t,\Omega)}\,dt\,.
$$

To show the converse inequality, we note that $(\mM_p(\Omega),Wb_p(\Omega))$ is a geodesic space. This is shown in \cite[Proposition 2.9]{FG10}) for $p=2$, minor modifications of the argument yield the result for general $p$. Let $(\mu_t)_{t\in[0,T]}$ be a geodesic curve with constant speed connecting $\mu_0$ and $\mu_T$. Then $(\mu_t)$ is absolutely continuous with $|\mu'|(t)=Wb_p(\mu_0,\mu_T)/T$. By Theorem \ref{teo1} there exists a time-dependent vector field $(\vv_t)$ such that $(\mu,\vv)\in \ce^\Omega_T$ and $\|\vv_t\|_{L^p(\mu_t)}\leq \frac{Wb_p(\mu_0,\mu_T)}T$ for a.e. $t.$
This immediately gives the desired inequality "$\geq$".
\end{proof}

\section{Non-linear diffusion equations with Dirichlet bc as gradient flows}\label{sec:non-lin-gf}

In this section, we consider gradient flows with respect to the distance $Wb_2$ of  internal energy functionals of the form
\[
\cF(\mu) = \int F\Big(\frac{d\mu}{d\Leb_\Omega}\Big)d \Leb_\Omega\;,
\]
for suitable strictly convex functions $F:[0,\infty)\to[0,\infty)$. We will show that the $Wb_2$-gradient flow (in the sense of curves of maximal slope) is given by $\mu_t=\rho_t\Leb_\Omega$ with $\rho_t$ a solution to the following non-linear Cauchy problem with Dirichlet boundary conditions $\lambda\ge0$,
\begin{equation}\label{eq:NLD}
\begin{cases}
\partial_t\rho = \Delta L_F(\rho) & \text{in }  \Omega\times(0,+\infty) \\
\rho(0,\cdot) = \rho_0 & \text{in } \Omega \\
\rho = \lambda & \text{on } \partial\Omega\times(0,\infty)\;.
\end{cases}
\end{equation}
 Here, $L_F:[0,\infty)\to[0,\infty)$ is the pressure function associated to $F$ defined by
\begin{equation}\label{LF}
L_F(r)= rF'(r) -F(r)\;,
\end{equation}
and $\lambda\geq 0$ is the unique point where $F$ is minimal. We will start by introducing and analysing the energy functional and its dissipation. Then we characterise curves of maximal slope as solutions to the above PDE.

\subsection{Internal energy and dissipation}~\label{InEnDis}
\smallskip

We fix an internal energy density $F$ satisfying the following properties.

\begin{assumption}\label{ass:F}
The function $F:[0,\infty)\to[0,\infty)$ satisfies 
\begin{enumerate}
\item[(i)] $F$ is continuous on $[0,\infty)$ and $C^2$ on $(0,\infty)$ with $F''>0$;
\item[(ii)] $F$ is superlinear, i.e.
\[
\lim_{r\to\infty} \frac{F(r)}{r} = \infty\;;
\]
\item[(iii)] $F$ satisfies the doubling condition, i.e. there exists $C>0$ such that 
\[F(r+s) \leq C \big(1+F(r)+F(s)\big)\quad \forall r,s\geq 0\;;\]
\item[(iv)] $r\mapsto \sqrt{r}F''(r)$ is locally integrable on $[0,\infty)$; the function $G(r):=\int_0^r\sqrt{s}F''(s)d s$ is such that there is a constant $C>0$ with 
\begin{equation}\label{ineq-G}
\sqrt{r}\leq C\big(1+G(r)\big)\quad \text{for all}\,\,\, r\ge0;
\end{equation}
\item[(v)] the function $h:(0,\infty)\to(0,\infty)$ defined as 
\begin{equation}\label{def-h}
h(r):=\big(\sqrt{r}F''(r)\big)^{-2}
\end{equation}
 is concave;
\item[(vi)] $F$ attains its (unique) minimum at $\lambda\geq 0$ and $F(\lambda)=0$;
\end{enumerate}

\end{assumption}

\begin{remark}\label{rem:mod-F}
If a function $F:[0,\infty)\to\R$ satisfying properties (i)-(v) 
above and $\lambda>0$ are given, we can construct a function $F_\lambda:[0,\infty)\to[0,\infty)$ satisfying additionally also property (vi) by setting 
\[F_\lambda(r)= F(r) - F(\lambda)-F'(\lambda)(r-\lambda)\;.\]
Indeed, properties (i)-(v), remain valid for $F_\lambda$ and due to the strict convexity of $F$, we have $F_\lambda\geq 0$ and $F_\lambda(r)=0$ if and only if $r=\lambda$, i.e. (vi) holds. Provided that $F$ is (right-) differentiable at $0$ we can proceed in the same way also if $\lambda=0$ is given. We further note that $L_{F_\lambda}(r)=L_F(r)-L_F(\lambda)$.
\end{remark}

\begin{remark}\label{limitLF}
Under the above assumptions on $F$, the pressure function $L_F$ given by \eqref{LF} has the following properties:
There is a constant $C>0$ such that 
\[0\leq L_F(r)\leq C(1+F(r))\qquad \forall r\geq 0\;.\]
Moreover, the function $L_F:[0,\infty)\to [0,\infty)$ is strictly increasing and unbounded.

Indeed, the convexity of $F$ and the doubling condition yield
\[0\leq rF'(r)-F(r)\leq F(2r)-2F(r)\leq C(1+F(r))\]
for all $r\geq 0$ and a suitable constant $C$. 
In particular, we have $F'(r)\geq F(r)/r$ and hence $F'(r)\nearrow +\infty$ as $r\nearrow+\infty$ by the superlinearity of $F$. Moreover, $L_F'(r)=rF''(r)>0$ for all $r$, so that $L_F$ is strictly increasing. To see that $L_F$ is unbounded, note that $F'$ is increasing and hence
\begin{align*}
L_F(r)= rF'(r)-F(r) = \int_0^r \left[F'(r)-F'(s)\right]ds \geq \int_0^1\left[F'(r)-F'(s)\right]ds = F'(r)-F(1)+F(0)\;.
\end{align*} 
Thus, $L_F(r)\nearrow +\infty$ as $r\nearrow +\infty$.
\end{remark}

\begin{definition}[Internal energy]\label{def:energy} The \emph{internal energy} $\cF$ is defined for $\mu\in \mathcal M_+(\Omega)$ by 
\begin{equation}\label{eq51}
\cF(\mu) = \int F(\rho)\, d\Leb_\Omega,
\end{equation}
provided that $\mu=\rho\Leb_\Omega$ is absolutely continuous w.r.t. $\Leb_\Omega$. Otherwise we set $\mathcal F(\mu)=+\infty$.
\end{definition}

\begin{remark}\label{rem52}
Typical choices of the functions $F$ satisfying properties (i)-(iv) in Assumption \ref{ass:F} include the following:
\begin{enumerate}
\item[i)] $F(r)=r^\alpha/(\alpha-1)$ with $\alpha>1$.
\noindent The pressure is given by $L_F(r)=r^\alpha$ leading to the porous medium equation $\partial_t \rho=\Delta \rho^\alpha$.
For $\lambda\geq0$, we find 
\[
F_\lambda(r)= \frac{1}{\alpha-1}\big[r^\alpha - r\alpha \lambda^{\alpha-1}\big] +\lambda^\alpha\;.
\]
\item[ii)] $F(r)=r\log r$.

\noindent This corresponds to the limit case $\alpha\to 1$. The pressure is given by $L_F(r)=r$, leading to linear diffusion $\partial_t\rho=\Delta\rho$ in \eqref{eq:NLD}.
For $\lambda>0$, we find 
\[
F_\lambda(r)=r\log(r/\lambda) -r+\lambda\;.
\]
As $F$ has infinite slope at 0, the case $\lambda=0$ cannot be considered.
\end{enumerate}

Property (v) is satisfied for ii) and for i) provided $\alpha \leq 3/2$.
\end{remark}

\begin{remark}\label{rem:F-lsc}
Note that due to the superlinearity of $F$, for $\mu$ with $\mathcal F(\mu)<\infty$ we have that $\mu=\rho\Leb_ \Omega $ with $\rho\in L^1(\Omega)$ and hence $\mu(\Omega)<\infty$. The functional $\cF:\mathcal M_+(\Omega)\to[0,+\infty]$ is lower semicontinuous with respect to vague convergence. This follows from the convexity and superlinearity of $F$ by general results on integral functionals, see e.g.~\cite[Cor. 3.4.2]{But89}. (Note that this reference considers functional defined on finite measure, while $\mathcal M_+(\Omega)$ is the space of locally finite measures. However, due to the previous comment, only sequences of finite measures need to be considered for the lower semicontinuity.)
\end{remark}

Let us define the function $G:[0,\infty)\to[0,\infty)$ by setting
\begin{equation}\label{funcG}
G(r)=\int_0^r \sqrt{s}F''(s)ds\;.
\end{equation}
Note that $G$ is strictly increasing and we have the identity 
\begin{equation}\label{Gprimo}
G'(r) = \sqrt{r} F''(r) = \frac{L_F'(r)}{\sqrt{r}}\;.
\end{equation}

\begin{definition}\label{endiss}
The \textnormal{energy dissipation functional} $\cI:\mathcal{M}_+(\Omega)\to[0,+\infty]$ is defined by
\begin{equation}\label{eq:def-diss}
\cI(\mu):=\begin{cases}
\int_{\Omega}|\nabla G(\rho)|^2\, d\Leb\;&\quad\text{if}\,\,\,\mu=\rho \Leb_\Omega\,\,\,\text{with}\,\,\, G(\rho)\in W^{1,2}(\Omega);\\
+\infty&\quad\text{otherwise}\end{cases}
\end{equation}
Moreover, we define the functional $\overline{\cI}:\mathcal{M}_+(\Omega)\to[0,+\infty]$ by setting
$$
\overline{\cI}(\mu) = \begin{cases} \cI(\mu) &\quad\text{if}\,\,\mu=\rho\Leb_\Omega\,\,\text{with}\,\,G(\rho)\in W^{1,2}(\Omega)\,\,\text{and}\,\,\mathcal T (G(\rho))=G(\lambda);\\
+\infty &\quad\text{otherwise}\,.\end{cases}$$
Namely, $\overline{\cI}(\mu)$ and $\cI(\mu)$ coincide if $G(\rho)-G(\lambda)\in W_0^{1,2}(\Omega)$.
\end{definition}
Recall the function $h(r)=\big(\sqrt{r}F''(r)\big)^{-2}$ defined in \eqref{def-h} (see Assumption \ref{ass:F} (v)).

\begin{lemma}\label{lem:reg-rho}
Let $\mu=\rho\Leb_\Omega$ with $\mathcal I(\mu)<\infty$, i.e.~$G(\rho)\in W^{1,2}(\Omega)$. Then we have $\rho\in W^{1,1}(\Omega)$ with $\nabla \rho =\sqrt{h(\rho)}\nabla G(\rho)$. Moreover, $\rho$ has trace $\lambda$ on $\partial \Omega$ if and only if $G(\rho)$ has trace $G(\lambda)$.
\end{lemma}

\begin{proof}
 As a direct consequence of \eqref{ineq-G} in (iv), Assumption \ref{ass:F} and $g:=G(\rho)\in L^2(\Omega)$ we have that $\rho\in L^1(\Omega)$.  Recall that $G$ is $C^1$ and strictly increasing in $(0,\infty)$, see \eqref{funcG} and \eqref{Gprimo}. Define its inverse function $H:=G^{-1}$ and approximate it by Lipschitz functions $H_n$ setting $$H_n(r) :=\int_0^r\min\{H'(s),n\}\,ds.$$ Then we deduce that $\rho_n:=H_n(g)\in W^{1,1}(\Omega)$ with 
$$\nabla \rho_n = H_n'(g)\nabla g = \min\{n,(G^{-1})'(g)\}\nabla g= \min\{n, \sqrt{h(\rho)}\}\nabla g,$$
with $h$ as in \eqref{def-h}. Note that $\rho_n\to \rho \in L^1(\Omega)$. From concavity of $h$, Assumption \ref{ass:F}-(v), we deduce that $h(\rho)\in L^1(\Omega)$. Hence,  the fact that $\nabla g\in L^2(\Omega)$ and Cauchy-Schwarz inequality together with dominated convergence theorem allow us to conclude that $\rho\in W^{1,1}(\Omega)$ with $\nabla \rho = \sqrt{h(\rho)}\nabla G(\rho)$. If $G(\rho)$ has trace $G(\lambda)$ one first checks by approximation with smooth functions that $\rho_n$ has trace $H_n(G(\lambda))$. Since $\rho_n\to\rho\in W^{1,1}(\Omega)$ the continuity of the trace operator gives that $\rho$ has trace $\lambda$. The reverse implication is treated similarly.
 \end{proof}

\begin{remark}\label{rem:GvsL}
Let $\mu=\rho\Leb_\Omega$ with $\mathcal F(\mu)<\infty$. Then $G(\rho)\in W^{1,2}(\Omega)$ if and only if  
$L_F(\rho)\in W^{1,1}(\Omega)$ with $\nabla L_F(\rho)=\bw \rho$ for a vector field $\bw\in L^2(\mu)$. In this case, we have
\[
\cI(\mu)=\int_{\Omega}|\bw|^2 \,d\mu\;.
\]
Indeed, recall from Remark \ref{rem:F-lsc} that the assumption $\cF(\mu)<\infty$ implies that $\rho\in L^{1}(\Omega)$. Further, from Remark \ref{limitLF} we infer that  $L_F(\rho)\in L^1(\Omega)$. Now, the claim can be checked starting from the observation that formally we have
\[
\nabla L_F(\rho)= \sqrt{\rho}\nabla G(\rho)\;.
\]
If $G(\rho)\in W^{1,2}(\Omega)$, we deduce that $\nabla L_F(\rho)$ in $L^1(\Omega)$ by Cauchy-Schwarz since $|\nabla G(\rho)|$ and $\sqrt{\rho}$ belong to $L^2(\Omega)$. Further $\nabla L_F(\rho)=\bw \rho$ with $\bw=\nabla G(\rho)/\sqrt{\rho} \in L^2(\mu)$. Conversely, if $L_F(\rho)\in W^{1,1}(\Omega)$ with $\nabla L_F(\rho)=\bw \rho$ and $\bw\in L^2(\mu)$, we deduce that $\nabla G(\rho) = \sqrt{\rho}\bw\in L^2(\Omega)$. By Poincar\'e inequality, this implies $G\in W^{1,2}(\Omega)$.

\end{remark}

Next, we show that the functionals $\cI$ and $\overline{\cI}$ are lower semicontinuous w.r.t. vague convergence on sublevels of $\cF$.

\begin{lemma}\label{lem:lsc-diss}
 Let $\{\mu_n\}\subset\mathcal{M}_+(\Omega)$ be a sequence of measures vaguely converging to $\mu\in\mathcal{M}_+(\Omega)$ such that
\begin{equation}\label{hplem2}
\sup_{n}\,\cF(\mu_n)<+\infty\;\quad\text{and}\quad \sup_{n}\,\cI(\mu_n)<\infty\;.
\end{equation}
Then we have
\begin{equation}\label{eq:I-lsc}
\cI(\mu)\le\liminf_{n}\cI(\mu_n)\;.
\end{equation}
If in addition $\overline{\cI}(\mu_n) <\infty$ for all $n$, then also $\overline{\cI}(\mu) <\infty$, i.e. \eqref{eq:I-lsc} holds with $\overline{\cI}$ instead of $\mathcal I$ as well.
\end{lemma}

\begin{proof}
For any $n\in\mathbb{N}$, let $\mu_n=\rho_n\Leb_\Omega$. Due to \eqref{hplem2} and the lower semicontinuity of $\mathcal{F}$ we have that $\cF(\mu)<\infty$.
Therefore we can write $\mu=\rho\Leb_\Omega$  for a suitable $\rho$. Superlinearity of $F$ implies that  
$\rho_n$ converges weakly to $\rho$ in $L^1(\Omega)$. 
%The doubling condition and convexity of $F$ imply that $0\leq L_F(r)\leq F(2r)\leq C\big(1+F(r)\big)$ and thus \eqref{hplem2} yields that $L_F(\rho_n)$ is bounded in $L^1(\Omega)$. We can extract a subsequence such that $L_F(\rho_n)\to L$ in $L^1_{\operatorname{loc}}(\Omega)$ and pointwise a.e. By truncation and monotonicity of $L_F$ we deduce that $L=L_F(\rho)$ a.e. and therefore $\nabla L_F(\rho_n)\to\nabla L_F(\rho)$ in the sense of distributions and $\rho_n\to\rho$ in $L^1(\Omega)$.
Due to \eqref{eq:def-diss}, \eqref{hplem2}, and Poincar\'e inequality we have
\[
\sup_{n}\|G(\rho_n)\|_{W^{1,2}(\Omega)}<\infty\;.
\]
By weak compactness and the Rellich-Kondrachov theorem there exists $g\in W^{1,2}(\Omega)$ such that, up to a subsequence, we have that $G(\rho_{n})\rightharpoonup g$ weakly in $W^{1,2}(\Omega)$, $G(\rho_n)\to g$ strongly in $L^2(\Omega)$ and pointwise a.e. From the monotonicity of $r\mapsto G(r)$ and the fact that $\rho_n\rightharpoonup \rho$ weakly in $L^1$, we deduce that $g=G(\rho)$. From the weak lower semicontinuity of the $L^2$-norm we get
\[
\liminf_{n\to\infty}\cI(\mu_n)=\liminf_{n\to\infty} \int_\Omega |\nabla G(\rho_n)|^2d\Leb \geq \int_\Omega|\nabla G(\rho)|^2d\Leb=\cI(\mu)\,.
\]
Thus \eqref{eq:I-lsc} is established.

Assume now that additionally that $\overline{\cI}(\mu_n) <\infty$. Then $G(\rho_n)-G(\lambda)\in W^{1,2}_0(\Omega)$,  i.e. $G(\rho_n)$ has trace $G(\lambda)$ on $\partial\Omega$ and $\overline{\cI}(\mu_n) =\mathcal I(\mu_n)$. We show that $G(\rho)$ has trace $G(\lambda)$ using the characterisation of the trace operator $\mathcal T$ via an integration by parts formula, see \cite[Theorem 4.6]{EvGa}. Since $G(\rho_n)$, $G(\rho)\in W^{1,2}(\Omega)$, we have for any smooth and bounded vector field $\Phi$ on $\bar\Omega$ that
\[
\int_\Omega G(\rho_n)\nabla\cdot\Phi\,d\Leb+\int_\Omega\nabla G(\rho_n)\cdot\Phi\,d\Leb=\int_{\partial\Omega}(\Phi\cdot \nu)\mathcal{T}\left(G(\rho_n)\right)\,d\mathcal{H}^{d-1}=\int_{\partial\Omega}(\Phi\cdot \nu)G(\lambda)\,d\mathcal{H}^{d-1},
\]
and
\[
\int_\Omega G(\rho)\nabla\cdot\Phi\,d\Leb+\int_\Omega\nabla G(\rho)\cdot\Phi\,d\Leb=\int_{\partial\Omega}(\Phi\cdot \nu)\mathcal{T}\left(G(\rho)\right)\,d\mathcal{H}^{d-1}\;.
\]
From the weak convergence $G(\rho_n)\rightharpoonup G(\rho)$ in $W^{1,2}(\Omega)$ we infer that
\[
\int_{\partial\Omega}(\Phi\cdot \nu)G(\lambda)\,d\mathcal{H}^{d-1}=\int_{\partial\Omega}(\Phi\cdot \nu)\mathcal{T}\left(G(\rho)\right)\,d\mathcal{H}^{d-1}\;,
\]
and hence, by arbitrariness of $\Phi$, we have $\mathcal{T}(G(\rho))=G(\lambda)$. Thus $G(\rho)-G(\lambda)\in W^{1,2}_0(\Omega)$ and $\overline{\cI}(\mu) =\cI(\mu)<\infty$.
 \end{proof}

\begin{lemma}\label{lem:I-cvx}
The functional $\mathcal I$ is convex w.r.t. linear interpolation, i.e. for any $\mu_0,\mu_1\in \mathcal M_+(\Omega)$ we have that, 
\begin{equation}\label{eq:I-cvx}
\cI(\mu_{t}) \leq (1-t)\,\cI(\mu_0) + t\,\cI(\mu_1)\;,
\end{equation}
for $\mu_{t}:=(1-t)\mu_0+t\mu_1$ for any $t\in[0,1]$. The same holds for $\overline{\cI}$.
\end{lemma}
\begin{proof}
We can assume that $\cI(\mu_0)$, $\cI(\mu_1)<\infty$ and hence
 $\mu_0=\rho_0\Leb_\Omega$ and $\mu_1=\rho_1\Leb_\Omega$ for densities $\rho_0,\rho_1$ with $G(\rho_0)$, $G(\rho_1)\in W^{1,2}(\Omega)$. 
 By Remark \ref{rem:GvsL} we can write
\[\cI(\mu) = \int_\Omega |G'(\rho)|^2|\nabla \rho|^2\,d\Leb = \int_{\{\rho>0\}} \frac{|\nabla\rho|^2}{h(\rho)}\,d\Leb=\int_\Omega\alpha_2\big(|\nabla \rho|,h(\rho)\big)\,d\Leb\;,\]
with $h(r)$ as in \eqref{def-h} %=\big(\sqrt{r}F''(r)\big)^{-2}$ 
and $\alpha_2$ be the convex function defined in \eqref{eq:alpha}. By Assumption \ref{ass:F}-(v), $h$ is concave and thus the function $(v,r)\mapsto\alpha\big(|v|,h(r)\big)$ is convex on $\R^d\times[0,\infty)$ and \eqref{eq:I-cvx} follows. 
We now consider $\overline{\cI}$. Similarly, we can assume $\overline{\cI}(\mu_0)$, $\overline{\cI}(\mu_1)<\infty$ which means that $\mu_0=\rho_0\Leb_\Omega$ and $\mu_1=\rho_1\Leb_\Omega$ for densities $\rho_0,\rho_1$ with $G(\rho_0)-G(\lambda)$, $G(\rho_1)-G(\lambda)\in W^{1,2}_0(\Omega)$. Then, from the previous argument, it suffices to show that for the linear combination $\rho_t=(1-t)\rho_0+t\rho_1$ it holds $\mathcal T\big(G(\rho_t)\big) =G(\lambda)$.  By Lemma \ref{lem:reg-rho} we have that $\rho_0,\rho_1\in W^{1,1}(\Omega)$ with $\mathcal T(\rho_0)=\mathcal T(\rho_1)=\lambda$ and hence also $\mathcal T(\rho_t)=\lambda$. Applying again Lemma \ref{lem:reg-rho} gives that $\mathcal{T}(G(\rho_t))=G(\lambda)$.
\end{proof}

\subsection{Variational characterisation}
We will now give a variational characterisation of the diffusion equation with Dirichlet boundary conditions in problem \eqref{eq:NLD}. Weak solutions will be characterised via an energy dissipation equality or in other words as minimisers of a De Giorgi functional associated to the internal energy $\cF$ and the distance $Wb_2$.
\smallskip

We adopt the following notion of weak solution.
\begin{definition}\label{def:weak-sol}
A curve of densities $(\rho_t)_{t\in[0,T]}$ is called a weak solution to the non-linear diffusion equation $\partial_t\rho = \Delta L_F(\rho)$ with Dirichlet boundary condition $\lambda$ and initial datum $\rho_0$ if $t\mapsto \rho_t$ is vaguely continuous on $[0,T]$, $t\mapsto L_F(\rho_t)$ belongs to $L^1\big([0,T];W^{1,1}_0(\Omega)\big)$ and for all $\varphi\in C^\infty_c(\Omega)$ and $0\leq s<t$ it holds
\begin{equation}\label{eq:weak-sol}
\int_\Omega\varphi (\rho_t-\rho_s)\,d\Leb = \int_s^t \int_\Omega \Delta \varphi\, L_F(\rho_r) \,d\Leb\,dr\;.
\end{equation}
\end{definition}

Recall that by the assumptions on $F$, see Assumption \ref{ass:F}, $\lambda$ is the unique zero of $L_F$, hence the requirement $L_F(\rho)\in W^{1,1}_0(\Omega)$ indeed encodes the boundary condition $\rho=\lambda$ on $\partial \Omega$.
\medskip

For a $2$-absolutely continuous curve $[0,T]\ni t\mapsto\mu_t$ in $(\mM_2(\Omega),Wb_2)$ with $\mathcal{F}(\mu_0)<+\infty$ we define the functional
\begin{equation}\label{eq:defDG}
\mathcal{L}_T(\mu):=\mathcal{F}(\mu_T)-\mathcal{F}(\mu_0)+\frac12 \int_0^T\left(|\mu_t'|^2+\overline{\mathcal{I}}(\mu_t)\right)\,dt\;.
\end{equation}

\begin{theorem}\label{teo2}
For any $2$-absolutely continuous curve $(\mu_t)_{t\in[0,T]}$ in $(\mathcal M_2(\Omega),Wb_2)$ with $\mathcal{F}(\mu_0)<+\infty$ we have that $\mathcal L_T(\mu)\geq 0$. Moreover, $\mathcal{L}_T(\mu_t)=0$ if and only if $\mu_t=\rho_t\Leb_\Omega$ such that  $t\mapsto G(\rho_t)-G(\lambda)$ belongs to $L^2\big([0,T],W^{1,2}_0(\Omega)\big)$ and $(\rho_t)$ is a weak solution to the non-linear diffusion equation 
\[\partial_t\rho = \Delta L_F(\rho)\;.\]
In this case, we have $|\mu'_t|^2=\overline{\cI}(\mu_t)$ for a.e. $t$ as well as 
the energy-dissipation identity
\begin{equation}\label{eq:energy-id}
\cF(\mu_T) - \cF(\mu_0) =  -\int_0^T\overline{\cI}(\mu_t)\,dt\;.
\end{equation}
\end{theorem}

%Note that the requirement that $G(\rho)-G(\lambda)\in W^{1,2}_0(\Omega)$ encodes the boundary condition $\rho=\lambda$ on $\partial\Omega$. 

\begin{remark}
In the language of gradient flows in metric spaces which is discussed in the next subsection, the previous theorem states that weak solutions to \eqref{eq:NLD} are precisely the curves of maximal slope of $\cF$ in the metric space $(\mathcal M_2(\Omega),Wb_2)$ w.r.t. the strong upper gradient $\sqrt{\overline{\cI}}$.
\end{remark}

To prove the latter theorem, the main step is to establish a chain rule for the energy $\mathcal{F}$ along absolutely continuous curves in $(\mathcal{M}_2(\Omega),Wb_2)$. 
%Let us define the function $\alpha:[0,+\infty)\times\R\to[0,+\infty)$ as
%\begin{equation}\label{alpha2}
%\alpha(s,u):=\begin{cases}
%\frac{|u|^2}{s}\quad&\text{if}\,\,s>0\\
%0\quad&\text{if}\,\,s=0\,\,\text{and}\,\,u=0\\
%+\infty\quad&\text{otherwise}\,.
%\end{cases}
%\end{equation}

\begin{proposition}[Chain rule]\label{propchainrule}
Let $(\mu_t)_{t\in[0,T]}$ be a $2$-absolutely continuous curve in $(\mathcal{M}_2(\Omega),Wb_2)$ with $\mathcal{F}(\mu_0)<+\infty$ and $(\vv_t)$ be a family of vector fields such that $(\mu,\vv)\in \mathcal{CE}^\Omega_T$  and 
\begin{equation}\label{eq53}
\int_0^T\|\vv_t\|^2_{L^2(\mu_t)}\,dt<\infty\;,\qquad \int_0^T\overline{\cI}(\mu_t)\,dt\,\,<+\infty\;.
\end{equation}
Let $\rho_t$ be such that $\mu_t=\rho_t\Leb_\Omega$ and let $\bw_t\in L^2(\mu_t)$ be the vector field defined by $\nabla L_F(\rho_t)= \rho_t\bw_t$ as in Remark \ref{rem:GvsL}. Then the function $t\mapsto\mathcal{F}(\mu_t)$ is absolutely continuous and
\begin{equation}\label{eq52}
\frac{d}{dt}\mathcal{F}(\mu_t)=\int_{\Omega}\left\langle\bw_t,\vv_t\right\rangle d\mu_t\qquad \text{ for a.e. } t\,.
\end{equation}
\end{proposition}

Before entering into the proof, let us motivate \eqref{eq52} by an informal computation and see in particular how the boundary condition enters. Assuming $(\mu,\vv)$ are sufficiently smooth and decaying, we calculate
\begin{align*}
\frac{d}{dt}\mathcal{F}(\mu_t) 
&= 
\int_\Omega F'(\rho_t)\partial_t\rho_t\; d\Leb
 = 
- \int_\Omega F'(\rho_t) \nabla\cdot(\rho_t \vv_t)\; d\Leb\\
&= 
\int_\Omega \langle\nabla F'(\rho_t),\vv_t\rangle \rho_t\; d\Leb
=\int_\Omega\langle \bw_t,\vv_t\rangle\; d\mu_t\;.
\end{align*}
Here, we have used the continuity equation in the second step, integration by parts in the third step and in the last step that $\bw_t=\nabla F'(\rho_t)$ because $\rho_t\bw_t=\nabla L_F(\rho_t)=\rho_t\nabla F'(\rho_t)$. Note that no boundary term occurs, since the assumption $\overline{\cI}(\mu_t)<\infty$ gives that $\rho_t=\lambda$ on $\partial \Omega$ and thus $F'(\rho_t)=F'(\lambda)=0$ on $\partial\Omega$ since $\lambda$ is the minimizer of $F$. In the following we will make the previous computation rigorous via a careful regularisation procedure.

\begin{proof}
We proceed in several steps. First, we regularise the curve $(\mu,\vv)$ and the energy $\mathcal F$, then we prove a chain rule for the regularised quantities, and finally we pass to the limit as the regularisation vanishes.\medskip

\emph{Step 1: Regularisation.}\smallskip

By assumption \eqref{eq53}, Definition \ref{endiss}, and the choice of $\vv$, we have that
\begin{equation}\label{eq54}
\int_0^T\int_{\Omega}|\bw_t|^2d\mu_t\,dt\,<\,\infty\;,\quad \int_0^T\int |\vv_t|^2d\mu_tdt<\infty\;.
\end{equation}
Let us set $J_t=\vv_t\mu_t=U_t\Leb_\Omega$ with $U_t=\vv_t\rho_t$. and perform the following regularisation on the pair $(\mu,J)\in \ce^\Omega_T$.

Firstly, we regularise in space as in the proof of Theorem \ref{teo1} (see Step 5): Let $\eta\in C_c^\infty(\R^d)$ be such that $\eta\geq0$, $\operatorname{supp}(\eta)\subset B_1(0)$, and $\|\eta\|_{L^1}=1$. For $\varepsilon>0$ define $\eta_\varepsilon$ by $\eta_\varepsilon(x)=\varepsilon^d\eta(x/\varepsilon)$ and set
\[\Omega_\varepsilon = \{x\in \Omega: {\sf dist}(x,\partial\Omega)>\varepsilon\}\;.\]
For each $t$ we define a pair $(\mu_t^\varepsilon,J_t^\varepsilon)$ on $\Omega_{2\varepsilon}$ by setting
\[\mu_t^\varepsilon := (\mu_t*\eta_\varepsilon)\vert_{\Omega_{2\varepsilon}}\;,
\quad J_t^\varepsilon := (J_t*\eta_\varepsilon)\vert_{\Omega_{2\varepsilon}}\;.\]
Note that only $\mu_t\vert_{\Omega_\varepsilon}$ and $J_t\vert_{\Omega_\varepsilon}$ contribute to the above integrals.
We recall that $(\mu^\varepsilon,J^\varepsilon)\in\ce_T^{\Omega_{2\varepsilon}}$, i.e. the regularised measures satisfy the continuity equation in the smaller domain $\Omega_{2\varepsilon}$. Let us write $J_t^\varepsilon=U^{\varepsilon}_t\Leb=\vv_t^\varepsilon\mu_t^\varepsilon$ and $\mu_t^\varepsilon=\rho_t^\varepsilon \Leb_{\Omega_{2\varepsilon}}$. 
%Also recall from \eqref{eq113} that by convexity of the Benamou-Brenier functional
%\begin{equation}
%\int_0^T\int_{\Omega_{2\varepsilon}}|\vv_t^{\varepsilon}|^2\,d\mu_t^{\varepsilon}dt\le\int_0^T\int_{\Omega}|\vv_t|^2\,d\mu_tdt\;.\label{eq:action-reg}
%\end{equation}
Lemma \ref{lem:I-cvx} yields that 
\begin{equation}\label{eq:I-reg}
\int_0^T\int_{\Omega_{2\varepsilon}}|\nabla G(\rho^\varepsilon_t)|^2\,d\Leb dt\leq \int_0^T\int_\Omega |\nabla G(\rho_t)|^2\,d\Leb dt = \int_0^T\cI(\mu_t) dt\;.
\end{equation}
Secondly, we regularise in time. Let $(\xi^{\delta})_{\delta>0}$ be a family of regular, even, kernels with support in the interval $[-\delta,\delta]$. Then for any $t\in[0,T]$ we define,
\begin{equation}\label{eq56}
\mu_t^{\varepsilon,\delta}:=\int_{-\delta}^\delta\xi^\delta(s)\mu^\varepsilon_{t-s}\,ds\;,
\quad J_t^{\varepsilon,\delta}:=\int_{-\delta}^\delta\xi^\delta(s)J^\varepsilon_{t-s}\,ds\;.
\end{equation} 
Here we assume that $\mu$ and $J$ have been extended with constant values on $[-\delta,0]$ and $[T,T+\delta]$. We write $\mu_t^{\varepsilon,\delta}:=\rho_t^{\varepsilon,\delta}\Leb_{\Omega_{2\varepsilon}}$ and $J^{\varepsilon,\delta}_t=U^{\varepsilon,\delta}_t\Leb_{\Omega_{2\varepsilon}}=\vv^{\varepsilon,\delta}_t\rho^{\varepsilon,\delta}_t\Leb_{\Omega_{2\varepsilon}}$.
Observe that $(\mu^{\varepsilon,\delta},J^{\varepsilon, \delta})\in \ce^{\Omega_{2\varepsilon}}_T$. 
%Again by convexity of the Benamou-Brenier functional and Lemma \ref{lem:I-cvx} we deduce that
%\begin{equation}\label{eq512}
%\begin{split}
%\int_0^T\int_{\Omega_{2\varepsilon}}|\nabla G(\rho^{\varepsilon,\delta}_t)|^2d\Leb_{\Omega_{2\varepsilon}} dt&\leq \int_0^T\cI(\mu_t) dt\;,\\
%\int_0^T\int_{\Omega_{2\varepsilon}}|\vv_t^{\varepsilon,\delta}|^2\,d\mu_t^{\varepsilon,\delta}dt&\le\int_0^T\int_{\Omega}|\vv_t|^p\,d\mu_tdt\;.
%\end{split}
%\end{equation}
\smallskip

As a third regularisation, we approximate $F$ by a family of convex functions $(F_\alpha):[0,\infty)\to\R$ for $\alpha>0$. To this end, we set $(F^ \alpha)':=\min\big\{\alpha, \max\{ F',-\alpha\}\big\}$ and 
\begin{equation}\label{Falpha}
F^{\alpha}(z):=\int_\lambda^{z}(F^{\alpha})'(s)\,ds\;,
\end{equation}
where $\lambda$ is the boundary condition in \eqref{eq:NLD}. Note that $F^\alpha$ is again convex, non-negative and attains its unique minimum $0$ at $\lambda$. Moreover, $F^ \alpha\nearrow F$ as $\alpha\nearrow \infty$. Then we define the regularized internal energy $\mathcal{F}^{\alpha}$ by 
\begin{equation}\label{eq59}
\mathcal{F}^{\alpha}(\mu):=\int_{\Omega}F^{\alpha}(\rho)\,d\Leb .%\color{magenta}+ \theta_\lambda \mu_s(\Omega)\;\normalcolor.
\end{equation}
provided $\mu=\rho\Leb_ \Omega $  and $\mathcal F^ \alpha (\mu)=+\infty$ else. Let us set for any $t\in[0,T]$, $\varepsilon,\delta, \alpha>0$
\begin{equation}\label{g}
g^{\varepsilon,\delta,\alpha}:=(F^{\alpha})'(\rho_t^{\varepsilon,\delta}).
\end{equation} 
By construction, $g^{\varepsilon,\delta,\alpha}$ is $C^1$ on $\Omega_{2\varepsilon}$. 

Finally, consider an increasing family of cut-off functions $\phi_\kappa\in C^\infty_c(\Omega)$ for $\kappa>0$ such that $0\leq \phi_{\kappa'}\leq \phi_\kappa\leq 1$ for all $\kappa'>\kappa$, $\phi_\kappa=1$ on $\Omega_{4\kappa}$, $\phi_\kappa=0$ on $\Omega\setminus\Omega_{2\kappa}$, and $|\nabla\phi_\kappa|\leq 1/\kappa$. We then define the functional $\mathcal F^{\alpha,\kappa}$ by 
\begin{equation}\label{eq:cutoff-energy}
\mathcal{F}^{\alpha,\kappa}(\mu):=\int_{\Omega}\phi_\kappa\cdot F^{\alpha}(\rho)\,d\Leb\; .
\end{equation}
provided $\mu=\rho\Leb$  and $\mathcal F^{\alpha,\kappa} (\mu)=+\infty$ else.
\medskip

\emph{Step 2: Regularised ($\delta\varepsilon\kappa\alpha\mathcal F$) chain rule.}\smallskip

In what follows, we write integrals over all $\Omega$ even though $\rho^\varepsilon$ is only defined on $\Omega_{2\varepsilon}$ with the understanding that due to the presence of $\phi_\kappa$ which is supported in $\Omega_{2\kappa}\subset\Omega_{2\varepsilon}$, there is no contribution from $\Omega\setminus\Omega_{2\varepsilon}$.

We claim that for $\delta,\varepsilon,\kappa,\alpha>0$ with $\varepsilon<\kappa$ and any $t\in [0,T]$ we have
\begin{equation}\label{eq510}
\mathcal{F}^ {\alpha,\kappa}(\mu_t^{\varepsilon,\delta})-\mathcal{F}^{\alpha,\kappa}(\mu_0^{\varepsilon,\delta})=\int_0^t\int_{\Omega} \left\langle\nabla \big(\phi_\kappa g_r^{\varepsilon,\delta, \alpha}\big),\vv^{\varepsilon,\delta}_r\right\rangle\, d\mu^{\varepsilon,\delta}_rdr\;,
\end{equation}
%where we set $g^{\varepsilon,\delta} = F'(\rho^{\varepsilon,\delta})$. 
Indeed, note that due to the smoothing, $(\rho^{\varepsilon,\delta}, \vv^{\varepsilon,\delta})$ satisfy the continuity equation 
\[\partial_t\rho^{\varepsilon,\delta}+\nabla\cdot(\rho^{\varepsilon,\delta}\vv^{\varepsilon,\delta})=0\]
in the classical sense on $\Omega_{2\varepsilon}$. Further $g^{\varepsilon,\delta,\alpha}=(F^{\alpha})'(\rho^{\varepsilon,\delta})$ is bounded $C^1$ on $\Omega_{2\varepsilon}$, $\phi_{\kappa}$ is smooth and compactly supported in $\Omega_{2\kappa}\subset \Omega_{2\varepsilon}$
%and away from $0$ 
and $\partial_r\rho^{\varepsilon,\delta}$ is uniformly bounded. Thus, we can differentiate under the integral and use the continuity equation to obtain
\begin{align*}\label{eq:chain-reg-pre}
\frac{d}{dr} \mathcal F^{\alpha,\kappa} (\mu^{\varepsilon,\delta}_r) 
&= \int_{\Omega} \phi_\kappa(F^\alpha)'(\rho^{\varepsilon,\delta}_r)\partial_r\rho^{\varepsilon,\delta}_rd\Leb
= \int_{\Omega}\left\langle\nabla \big(\phi_\kappa g_r^{\varepsilon,\delta,\alpha}\big),\vv^{\varepsilon,\delta}_r\right\rangle \rho^{\varepsilon,\delta}_r\, d\Leb\;.
\end{align*}
Note that no boundary term occurs in the integration by parts in the last step, since $\phi_k$ is compactly supported in $\Omega_{2\varepsilon}$. Integrating in time on $[0,t]$ yields \eqref{eq510}.
\medskip

\emph{Step 3: Passing to the limit.}\smallskip

We now pass to the limit as $\delta\to 0$, $\varepsilon\to 0$, $\kappa\to0$, and $\alpha\to \infty$ in this order in the regularised chain rule \eqref{eq510}. 
\medskip

\emph{Limit $\delta\to0$:}\smallskip

Note that $\rho^{\eps,\delta}$, $\vv^{\varepsilon,\delta}$, $\nabla\big(\phi_\kappa g^{\eps,\delta,\alpha}\big)$, as well as $F^\alpha(\rho^{\varepsilon,\delta})$ are all uniformly bounded on $\Omega_{2\varepsilon}\times[0,T]$ as $\delta\to0$ and converge pointwise to $\rho^{\varepsilon}$, $\vv^{\eps}$, $\nabla\big(\phi_\kappa g^{\varepsilon,\alpha}\big)$, and $F^{\alpha}(\rho^\varepsilon)$ respectively with $g^{\eps,\alpha}:=(F^\alpha)'(\rho^{\eps})$. Hence we can pass to the limit in both sides of \eqref{eq510} by dominated convergence and obtain
\begin{equation}\label{eq:chain-eps}
\mathcal{F}^{\alpha,\kappa}(\mu_t^{\varepsilon})-\mathcal{F}^{\alpha,\kappa}(\mu_0^{\varepsilon})=\int_0^t\int_{\Omega} \left\langle\nabla \big(\phi_\kappa g_r^{\varepsilon,\alpha}\big),\vv^{\varepsilon}_r\right\rangle\, d\mu^{\varepsilon}_rdr\;.
\end{equation}

\emph{Limit $\eps\to0$:}\smallskip

We first consider the right hand side of \eqref{eq:chain-eps} and claim that as $\varepsilon\to 0$
\begin{equation}\label{eq:RHSeps}
\int_0^t\int_{\Omega} \left\langle\nabla \phi_\kappa \, g_r^{\varepsilon,\alpha}+\phi_\kappa\,\nabla g_r^{\eps,\alpha},\vv^{\varepsilon}_r\right\rangle\, d\mu^{\varepsilon}_rdr \longrightarrow \int_0^t\int_{\Omega} \left\langle\nabla\phi_\kappa\, g_r^\alpha+\phi_\kappa\nabla g_r^\alpha,\vv_r\right\rangle\, d\mu_r dr \;,
\end{equation}
with $g_r^\alpha=(F^\alpha)'(\rho_r)$. Recall that $\vv^\eps\rho^\eps=\boldsymbol{1}_{\Omega_{2\eps}}\eta_\eps*(\vv\rho)$ and thus $\rho^\eps$ and $\vv^\eps\rho^\eps$ converge to $\rho$ and $\vv\rho$ a.e. on $\Omega_{2\kappa}\times[0,T]$ as $\eps\to0$. Also, by Lemma \ref{lem:reg-rho} we have $\rho_r\in W^{1,1}(\Omega)$ for a.e. $r$. Thus $\nabla \rho^\eps$ converges to $\nabla\rho$ and hence also $\nabla g^{\eps,\alpha}=(F^\alpha)''(\rho^\eps)\nabla\rho^\eps$ converges to $(F^\alpha)''(\rho)\nabla\rho = \nabla g^\alpha$ a.e.~in $\Omega_{2\kappa}\times[0,T]$.
We consider the first summand in \eqref{eq:RHSeps}. We have the majorant
\begin{align*}
\big|\langle\nabla\phi_\kappa g^{\varepsilon,\alpha}_r,\vv_r^\eps\rangle\rho_r^\eps\big|&\leq \boldsymbol{1}_{\Omega_{2\kappa}}\frac{\alpha}{\kappa}|\vv_r^\eps|\rho_r^\eps \;.
\end{align*}
From the bound $2|\vv|\rho\leq |\vv|^2\,\rho+\rho$, together with \eqref{eq54} we infer that $|\vv|\rho$ is in $L^1([0,T]\times\Omega_{2\kappa})$ and hence $|\vv_r^\eps|\rho_r^\eps\to|\vv_r|\rho_r$ in $L^1([0,T]\times\Omega_{2\kappa})$. Thus, dominated convergence yields the convergence of the first summand on the left of \eqref{eq:RHSeps} to the respective one on the right. To show convergence of the second summand, observe that we have the majorant
\begin{align*}
 \Big|\left\langle\phi_\kappa\nabla g_r^{\varepsilon,\alpha},\vv^{\varepsilon}_r\right\rangle\rho^\eps_r\Big|
&\leq
 \boldsymbol{1}_{\Omega_{2\kappa}} \frac12\left[|\nabla g^{\eps,\alpha}_r|^2\rho^\eps_r + |\vv^\eps_r|^2\rho^\eps_r \right]
\leq
 \boldsymbol{1}_{\Omega_{2\kappa}} \left[\frac{|\nabla \rho^\eps_r|^2}{2h(\rho^\eps_r)} + \frac{|U^\eps_r|^2}{2\rho^\eps_r} \right]
=: \Phi_1^\eps + \Phi_2^\eps\;,
\end{align*}
where we have used that $|\nabla g^{\eps,\alpha}_r|^2=|\nabla\rho_r^\eps|^2(F''(\rho_r^\eps))^2$ and $h$ is defined in \eqref{def-h}.
Note that, as $\eps\to0$
\[
\Phi^\eps_1 \to \Phi_1:= \frac{|\nabla\rho|^2}{2h(\rho)}\;,\qquad \Phi_2^\eps\to \Phi_2:=\frac{|U|^2}{2\rho }\quad \text{a.e. on } \Omega_{2\kappa}\times[0,T]\;.
\]
To conclude it, it suffices, by the (extended) dominated convergence theorem, (see e.g. \cite[Chap.~4,
Thm.~17]{Roy88}), to show that 
\begin{equation}\label{eq:RHSeps2}
\int_0^t\int_\Omega \Phi^\eps_i \,d\Leb \,dr \longrightarrow \int_0^t\int_\Omega \Phi_i \,d\Leb \,dr\quad \text{ as } \eps \to 0\;,
\quad i=1,2\;.
\end{equation}
To see this, we note that the function $(u,s)\mapsto |u|^2/h(s)$ is convex on $\R^n\times(0,\infty)$ when $h$ is concave and recall that $\rho^\eps=\eta_\eps*\rho$, $U^\eps=\eta_\eps*U$, to obtain the estimates
\begin{align*}
\Phi_i^\eps \leq \boldsymbol{1}_{\Omega_{2\eps}} \eta_\eps* \Phi_i\;,\quad i=1,2\;.
\end{align*}
Now, \eqref{eq:RHSeps2} follows in turn by the dominated convergence theorem and the fact that $\Phi_i\in L^1([0,T]\times\Omega)$. Hence, we conclude the convergence \eqref{eq:RHSeps}.
\smallskip

Now consider the left hand side of \eqref{eq:chain-eps}. Definition \eqref{Falpha} of $F^{\alpha}$ and the construction of $\phi_\kappa$ allows us to write 
$$
\left|\phi_\kappa F^\alpha(\rho)\right|\leq\left|\int_\lambda^\rho(F^\alpha)'(s)\,ds\right|\le \alpha(\rho-\lambda)\,.
$$
With the latter bound and the fact that $\rho^\eps=\boldsymbol{1}_{\Omega_{2\eps}}(\eta_\eps*\rho)$ converges to $\rho$ a.e. and in $L^1$, we get, by dominated convergence theorem, that as $\eps\to0$
\begin{equation*}
\cF^{\alpha,\kappa}(\mu^\eps)
 = \int_{\Omega}\phi_\kappa F^\alpha(\eta_\eps*\rho)\,d\Leb \rightarrow \int_\Omega \phi_\kappa F^\alpha(\rho)\,d\Leb = \cF^{\alpha,\kappa}(\mu)\;.
\end{equation*}
Thus, we obtain
\begin{equation}\label{eq:chain-kappa}
\mathcal{F}^{\alpha,\kappa}(\mu_t)-\mathcal{F}^{\alpha,\kappa}(\mu_0)= \int_0^t\int_{\Omega} \big\langle\nabla\phi_\kappa g^\alpha_r,\vv_r\big\rangle\, d\mu_r dr+ \int_0^t\int_{\Omega} \phi_\kappa\big\langle\nabla g^\alpha_r,\vv_r\big\rangle\, d\mu_r dr\;.
\end{equation}

\emph{Limit $\kappa\to0$:}\smallskip

We claim that as $\kappa\to0$, \eqref{eq:chain-kappa} converges to
\begin{equation}\label{eq:chain-alpha}
\mathcal{F}^{\alpha}(\mu_t)-\mathcal{F}^{\alpha}(\mu_0)= \int_0^t\int_{\Omega}\big\langle\nabla g^\alpha,\vv_r\big\rangle\, d\mu_r dr\;.
\end{equation}

By monotone convergence (recall $F^\alpha\geq0$), we have $\cF^{\alpha,\kappa}(\mu)\nearrow\cF^\alpha(\mu)$ as $\kappa\searrow 0$ for any $\mu$. From the bound $|\phi_\kappa\langle\nabla g^\alpha,\vv\rangle\rho|\leq |\bw||\vv|\rho$ and \eqref{eq54}, by dominated convergence we can pass to the limit in the second term on the right of \eqref{eq:chain-kappa}. It remains to show that the first term on the right of \eqref{eq:chain-kappa} goes to $0$. 
The latter is a boundary term and its vanishing will be due to the fact that $F'(\lambda)=0$ and $\rho=\lambda$ on $\partial\Omega$. To make this rigorous, note that the support of $\nabla\phi_\kappa$ is contained in the region  $B_{4\kappa}:=\{x\in \Omega: {\rm dist}(x,\partial \Omega)<4\kappa\}$. Thus, by assumption on $\phi_\kappa$ and the Cauchy-Schwarz inequality, we have the estimate
\begin{align*}
\left|\int_{\Omega} \big\langle\nabla\phi_\kappa g^\alpha_r,\vv_r\big\rangle\, d\mu \right|
&\leq
\frac{1}{\kappa}\int_{B_{4\kappa}} (F^\alpha)'(\rho_r)|\vv_r|\rho_r\; d\Leb\\
&\leq
\left(\frac{1}{\kappa^2}\int_{B_{4\kappa}} \left|(F^\alpha)'(\rho_r)\sqrt{\rho_r}\right|^2 \; d\Leb\right)^{1/2}\left(\int_{B_{4\kappa}}|\vv_r|^2\rho_r\; d\Leb\right)^{1/2}\;. 
\end{align*}
Using \eqref{Falpha}, \eqref{ineq-G} in Assumption \ref{ass:F}-(iv) and Lemma \ref{lem:bdry-est}, we bound the first factor as
\begin{equation*}
\int_{B_{4\kappa}} \big|(F^\alpha)'(\rho_r)\sqrt{\rho_r}|^2 \; d\Leb \leq \alpha^2C \int_{B_{4\kappa}}| 1+ G(\rho)|^2 d\Leb
\leq \alpha^2 C\kappa^2\int_{B_{4\kappa}}|\nabla G(\rho)|^2 d\Leb\;.
\end{equation*}
Collecting the previous estimates and a further Cauchy-Schwarz inequality in the time integration, we obtain
\begin{align*}
\left | \int_0^t\int_{\Omega} \left\langle\nabla\phi_\kappa g^\alpha_r,\vv_r\right\rangle\, d\mu_r dr\right |
&\leq 
C\alpha^2\left(\int_0^t\int_{B_{4\kappa}} |\nabla G(\rho_r)|^2 \; d\Leb \,dr \right)^{1/2}\left(\int_0^t\int_{B_{4\kappa}}|\vv_r|^2\rho_r\; d\Leb\,dr\right)^{1/2}\;. 
\end{align*}
Recalling that $\int_\Omega |\nabla G(\rho_r)|^2 d\Leb=\int_\Omega |\bw_r|^2d\mu_r=\overline{\mathcal I}(\mu_r)$, the finiteness of the integrals \eqref{eq54} readily gives that both factors in the last expression vanish as $\kappa\to0$ since $B_{4\kappa}\searrow \emptyset$.\medskip
%\newpage

\emph{Limit $\alpha\to\infty$:}\smallskip

We claim that as $\alpha\to\infty$, \eqref{eq:chain-alpha} converges to
\begin{equation}\label{eq:chain}
\mathcal{F}(\mu_t)-\mathcal{F}(\mu_0)= \int_0^t\int_{\Omega}\big\langle\bw_r,\vv_r\big\rangle\, d\mu_r dr\;.
\end{equation}

Observing that $\nabla g^\alpha=(F^\alpha)''\nabla\rho = \boldsymbol{1}_{\{F'(\rho)\in[-\alpha,\alpha]\}}\nabla F'(\rho)=\boldsymbol{1}_{\{F'(\rho)\in[-\alpha,\alpha]\}}\bw$ and recalling that the integral on the right of \eqref{eq:chain} is finite by the assumption \eqref{eq54}, we readily pass to the limit in the right hand side of \eqref{eq:chain-alpha}. Since $F^\alpha\nearrow F$ for $\alpha\nearrow\infty$ we pass to the limit in the left hand side by monotone convergence.

%By convexity of the function $F^\alpha$, and recalling that $0\leq F^\alpha\leq F$ and $\supp \phi_\kappa\subset \Omega_{2\eps}$, we deduce that for all $t$
%\begin{align*}
%\cF^{\alpha,\kappa}(\mu^\eps_t)
%&=
%\int_{\Omega}\phi_\kappa F^\alpha(\eta_\eps*\rho_t)d\Leb
%\leq
%\int_{\Omega_{2\eps}}\eta_\eps*F(\rho_t)d\Leb 
%\leq \int_\Omega F(\rho)d\Leb
%= \cF(\mu_t)\;.
%\end{align*}
%Combining this with the weak convergence of $\mu^\eps$ to $\mu$ and the lower semicontinuity of $\cF$, yields that $\cF(\mu^\eps_t)\to\cF(\mu_t)$ as $\eps\to 0$ for all $t\in[0,T]$.

Finally,  since $\cF(\mu_0)$ is finite and the right hand side of \eqref{eq:chain} is finite, we can conclude that $\cF(\mu_t)$ is finite for all $t\in[0,T]$ and we have
\[
\mathcal{F}(\mu_t)-\mathcal{F}(\mu_0)=\int_0^t\int_{\Omega} \left\langle\nabla \bw_r,\vv_r\right\rangle\, d\mu_rdr\;.
\]
Repeating the argument now with $0$ replaced by $s\in [0,T]$ yields that for all $s,t\in[0,T]$ we have 
\[
\mathcal{F}(\mu_t)-\mathcal{F}(\mu_s)=\int_s^t\int_{\Omega} \left\langle\nabla \bw_r,\vv_r\right\rangle\, d\mu_rdr\;.
\]
Hence, the map $t\mapsto \cF(\mu_t)$ is absolutely continuous and \eqref{eq52} holds.

\end{proof}

With the chain rule for the internal energy at hand we can now give the proof of Theorem \ref{teo2}.

\begin{proof}[Proof of Theorem \ref{teo2}]
Let $(\mu_t)_{t\in[0,T]}$ be a $2$-absolutely continuous curve in $(\mathcal M_2(\Omega), Wb_2)$, let $(\vv_t)$ be its optimal velocity vector field and set $J_t=\vv_t\mu_t$ such that $(\mu,J)\in\ce$. We can assume that \eqref{eq53} holds as otherwise $\mathcal L_T(\mu)=+\infty$, see \eqref{eq:defDG}. From Proposition \ref{propchainrule}, Cauchy-Schwarz and Young inequalities we infer that 
\begin{align*}
\mathcal{F}(\mu_T)-\mathcal F(\mu_0)&=\int_0^T\int_{\Omega}\left\langle\bw_t,\vv_t\right\rangle \,d\mu_t dt\\
&\geq -\frac12\int_0^T\int_\Omega |\bw_t|^2\,d\mu_tdt - \frac12 \int_0^T\int_\Omega|\vv_t|^2\,d \mu_t dt\\
&= -\frac12\int_0^T\left(|\mu'_t|^2 + \overline{\cI}(\mu_t)\right) dt\;.
\end{align*}
Here we have used Remark \ref{rem:GvsL} and Theorem \ref{teo1} in the last step. This shows that $\mathcal L_T(\mu)\geq0$.

We turn to the second part of the statement. Assume that $\mathcal{L}_T(\mu_t)=0$. Then, in particular \eqref{eq53} holds, i.e. $\mu_t=\rho_t\Leb_\Omega$ and $t\mapsto G(\rho_t)-G(\lambda)$ belongs to $L^2([0,T];W^{1,2}_0(\Omega))$. Moreover, equality must hold in the inequalities above. This implies that $\vv_t=-\bw_t$ a.e. w.r.t. the measure $d\mu_t dt$. Hence we have $J_t=-\bw_t\mu_t =-\nabla L_F(\rho_t)\Leb$ for a.e. $t$. Now the continuity equation for $(\mu,J)$ shows that $(\rho)$ is a weak solution to $\partial_t\rho=\Delta L_F(\rho)$. 

Conversely, if $\mu_t=\rho_t\Leb_\Omega$ such $t\mapsto G(\rho_t)-G(\lambda)$ belongs to $L^2([0,T];W^{1,2}_0(\Omega))$ is a weak solution, we have that \eqref{eq53} holds and $(\mu,J)$ satisfies the continuity equation with $J_t:=-\nabla L_F(\rho_t)\Leb_\Omega=\bw_t\mu_t$. Theorem \ref{teo1} yields that $(\mu_t)$ is $2$-absolutely continuous with
\[
|\mu'_t| \leq \|\bw_t\|_{L^2(\mu_t)}=\sqrt{\overline{\cI}(\mu_t)}\;.
\]
 From the chain rule (Proposition \ref{propchainrule}) together with the first part of the proof, we infer that 
\[ -\int_0^T\|\bw_t\|^2_{L^2(\mu_t)} \,dt = \cF(\mu_T)-\cF(\mu_0) \geq -\frac12 \int_0^T\|\bw_t\|^2_{L^2(\mu_t)}\,dt -\frac12 \int_0^T |\mu'_t|^2\,dt\;.\]
To combine the last two estimates shows that, in fact, $|\mu'_t|=\|\bw_t\|_{L^2(\mu_t)}$ for a.e. $t$ and hence $\mathcal  L_T(\mu) =0$.
The energy identity \eqref{eq:energy-id} now follows immediately.
\end{proof}

\subsection{Gradient flow characterisation}\label{sec:gradflow}

In this section we give a characterisation of the non-linear diffusion equation with Dirichlet boundary conditions \eqref{eq:NLD} as a metric gradient flow of $\cF$ in the space $(\mathcal M_2(\Omega),Wb_2)$. First, we briefly recall basic notions regarding gradient flows in metric spaces.
\smallskip

Let $(X,d)$ be a complete metric space and let
$E:X\to(-\infty,\infty]$ be a function with proper domain, i.e.~the
set $D(E):=\{x:E(x)<\infty\}$ is non-empty. Recall from the beginning of Section \ref{sec3} the definition of absolutely continuous curves. The following notion plays the role of the modulus of the gradient in a metric setting.

\begin{definition}[Strong upper gradient]\label{def:upper-grad}
  A function $g:X\to[0,\infty]$ is called a \emph{strong upper
    gradient} of $E$ if for any $x\in AC\big((a,b);(X,d)\big)$ the
  function $g\circ x$ is Borel and
  \begin{align*}
 |E(x_s)-E(x_t)|~\leq~\int_s^tg(x_r)|x'_r| d  r \quad\forall~a\leq s\leq t\leq
  b\ .
 \end{align*}
\end{definition}

Note that by the definition of strong upper gradient and Young's
inequality ($ab\leq \frac12(a^2+b^2)$), we have that for all $s\leq t$:
\begin{align*}
    E(x_t) - E(x_s) +\frac12 \int_s^t \left(g(x_r)^2 + |x'_r|^2\right)\, d  r\geq 0\;.
\end{align*}

\begin{definition}[Curve of maximal slope]\label{def:curve-max-slope}
  A locally $2$-absolutely continuous curve $(x_t)_{t\in(0,\infty)}$ is
  called a curve of maximal slope of $E$ w.r.t.~its strong upper
  gradient $g$, if $t\mapsto E(x_t)$ is non-increasing and
  \begin{align}\label{eq:cms}
    E(x_t) - E(x_s) +\frac12 \int_s^t \left(g(x_r)^2 + |x'_r|^2\right)\, d  r \leq 0 \quad\forall~0< s\leq t\;.
  \end{align}
  We say that a curve of maximal slope starts from $x_0\in X$ if $\lim_{t\searrow 0}x_t=x_0$.
\end{definition}

Equivalently, we can require equality in \eqref{eq:cms}. If a strong
upper gradient $g$ of $E$ is fixed, we also call a curve of maximal
slope of $E$ (relative to $g$) a \emph{gradient flow curve}.

Finally, we recall the definition of the (descending) local metric slope of $E$ as the
function $|\partial E|:D(E)\to[0,\infty]$ given by
\begin{align}\label{eq:metric-slope-def}
  |\partial E| (x) = \limsup_{y\to x}\frac{\max\{E(x)-E(y),0\}}{d(x,y)}\;.
\end{align}
The local metric slope is in general only a weak upper gradient $E$, see \cite[Thm.~1.2.5]{AGS}. We will also need the notion of relaxed slope $|\partial^-E|$, which is the sequentially $\sigma$-lower semicontinuous relaxation of the slope $|\partial E|$ w.r.t. a topology $\sigma$ on $E$ that could be weaker than convergence in the distance $d$. More precisely, one sets
\begin{equation}\label{eq:relslope}
|\partial^-E|(x):=\inf\left\{\liminf_n |\partial E|(x_n)\,\,:\,\,x_n\overset{\sigma}{\longrightarrow} x\right\}\,.
\end{equation}

\medskip

Let us now consider these notions in the metric space $(\mathcal M_2(\Omega),Wb_2)$ where the topology $\sigma$ is the topology of vague convergence. 

\begin{corollary}\label{cor}
$\sqrt{\overline{\cI}}$ is a strong upper gradient for $\mathcal{F}$ on $(\mathcal{M}_2(\Omega),Wb_2)$. A curve $(\mu_t)_{t\in[0,T]}$ is a curve of maximal slope w.r.t. this strong upper gradient if and only if $\mu_t=\rho_t\Leb_\Omega$ with $t\mapsto G(\rho_t)-G(\lambda)$ belonging to $L^2([0,T];W^{1,2}_0(\Omega))$ such that $(\rho_t)$ is a weak solution to \eqref{eq:NLD}.
\end{corollary}

\begin{proof}
For a $2$-absolutely continuous curve $[0,T]\ni t\mapsto \mu_t\in \mathcal{M}_2(\Omega)$ with optimal velocity vector field $(\vv_t)$ satisfying \eqref{eq53} we may apply Proposition \ref{propchainrule} and Cauchy-Schwarz inequality (as in the proof of Theorem \ref{teo2}) to get for any $0\leq s\leq t\leq T$:
\begin{align*}
|\mathcal{F}(\mu_t)-\mathcal{F}(\mu_s)|
&=
\Big|\int_s^t\int_\Omega\langle\bw_r,\vv_r\rangle d\mu_r dr\Big|
\le
\int_s^t\|\bw_r\|_{L^2(\mu_r)}\|\vv_r\|_{L^2(\mu_r)}dr\\
&=
\int_s^t\sqrt{\overline{\mathcal{I}}(\mu_r)}\,|\mu_r'|\,dr\,;
\end{align*}
where, in the last step, we have used Remark \ref{rem:GvsL} and Theorem \ref{teo1}. The latter shows that $\sqrt{\overline{\cI}}$ is a strong upper gradient. A curve $(\mu_t)$ is a curve of maximal slope w.r.t. this strong upper gradient if and only if $\mathcal L_T(\mu)=0$. Hence, the characterisation of curves of maximal slope follows by Theorem \ref{teo2}.
\end{proof}

In the following proposition we relate $\sqrt{\overline{\cI}}$ to the local slope of $\mathcal{F}$. This result shows that the Dirichlet boundary condition for the non-linear diffusion equation are encoded in the metric properties of $\cF$ w.r.t. the distance $Wb_2$.

\begin{proposition}\label{prop4}
For any  $\mu\in\mathcal{M}_2(\Omega)$ we have the bound
\begin{equation}\label{eq:slope-Fisher}
|\partial \mathcal{F}|^2(\mu)\ge \overline{\cI}(\mu)\,.
\end{equation}
\end{proposition}

\begin{remark}
We conjecture that equality holds in \eqref{eq:slope-Fisher}.
We recall that the local slope at $\mu$ of the functional $\cF$ w.r.t. the Wasserstein distance is well known to be equal to $\sqrt {\cI(\mu)}$, provided that the domain $\Omega$ is convex and $F$ satisfies the so called McCann conditions, see  \cite[Theorem 10.4.9]{AGS}. This crucially relies on the fact, that under these assumptions, $\cF$ is convex along Wasserstein geodesics. A serious source of difficulty in the present situation is the lack of convexity of $\cF$ (even in the case of the Boltzmann entropy) along geodesics in the modified Wasserstein distance $Wb_2$, see \cite[Rem.~3.4]{FG10}. The crucial feature of \eqref{eq:slope-Fisher} is that it holds with $\overline{\cI}$ in the right hand side, that is, finiteness of the local slope $|\partial \cF|(\mu)$ implies that $\mu=\rho\Leb_\Omega$ with $\rho$ satisfying the Dirichlet boundary condition $\rho=\lambda$ on $\partial\Omega$. This means that the boundary condition for the metric gradient flow is a consequence of the interaction of internal energy and the transport geometry. 
Let us point out that a similar result has been obtained using similar arguments in \cite[Thm.~1]{CMS} in the context of sticky-reflecting diffusions.
\end{remark}

\begin{proof}
\emph{Step 1:} We will first show the lower bound
\begin{equation}\label{eq:slope-Fisher-pre}
|\partial \mathcal{F}|^2(\mu)\ge \cI(\mu)\qquad \forall \mu\in \mathcal {M}_2(\Omega)\;.
\end{equation}
Fix $\mu\in \mathcal M_2(\Omega)$. We can assume that $\mu=\rho\Leb$ with $\rho\in L^1(\Omega)$ since otherwise $\cF(\mu)=+\infty$ and thus $|\partial\cF|(\mu)=+\infty$. We recall the following representation of the local slope, see \cite[Lem.~3.1.5]{AGS}:
\begin{equation}\label{eq:dual-slope}
\frac{1}{2}|\partial \cF|^2(\mu) = \limsup_{t\searrow 0} \frac{1}{t}\sup_{\nu}\Big[\cF(\mu)-\cF(\nu)-\frac{1}{2t}Wb_2(\mu,\nu)^2\Big]\;. 
\end{equation}
 By Lemma \ref{lem:weaksol-energy} there exists a weak solution $(\rho_t)$ to the non-linear diffusion equation \eqref{eq:NLD} with initial datum $\rho_0=\rho$ such that $t\mapsto G(\rho_t)-G(\lambda)$ belongs to $L^2\big([0,T],W^{1,2}_0(\Omega)\big)$ for any $T>0$. From Theorem \ref{teo2} we have that $|\mu'_s|^2=\cI(\mu_s)$ for a.e. $s$ and
 
\[
\cF(\mu_t)-\cF(\mu) = -\int_0^t \cI(\mu_s)\,ds\;.
\]
From Theorem \ref{teo1} and Remark \ref{rem1} we obtain that 
\[
Wb_2(\mu,\mu_t)^2 \leq t \int_0^t \cI(\mu_s)\,ds\;.
\]
Choosing $\nu=\mu_t$ in \eqref{eq:dual-slope} and combining it with the last two observations we obtain
\begin{align*}
\frac{1}{2}|\partial \cF|^2(\mu) 
&\geq
\limsup_{t\searrow 0} \frac{1}{t}\left[\cF(\mu)-\cF(\mu_t)-\frac{1}{2t}Wb_2(\mu,\mu_t)^2\right]\\
&\ge \limsup_{t\searrow 0} \frac{1}{2t}\int_0^t\cI(\mu_s)\,ds= \limsup_{t\searrow 0}\frac12 \int_0^1\cI(\rho_{ts})\, ds \geq \frac{1}{2}\,\cI(\mu)\;. 
\end{align*}
In the last step we have used Fatou's lemma together with the fact that $\mu_t$ converges vaguely to $\mu$ as $t\to 0$ and that $\cI$ is lower semi-continuous w.r.t. vague convergence. Note that the latter bounds hold regardless of whether $\cI(\mu)$ is finite or not.
\medskip

\emph{Step 2:} We will now show that the slope of $\cF$ at $\mu=\rho\Leb_\Omega$ is infinite unless $\rho=\lambda$ on $\partial\Omega$. More precisely, consider $\mu=\rho\Leb_\Omega$ with $G(\rho)\in W^{1,2}(\Omega)$ but $G(\rho)-G(\lambda)\notin W^{1,2}_0(\Omega)$, i.e. the trace of $G(\rho)$ on $\partial \Omega$ is not equal to the constant function $G(\lambda)$. Then we show that $|\partial \mathcal F|(\mu)=+\infty$.\smallskip

To this end, we consider the following perturbation of $\mu$: we define the measure 
\[\mu_\varepsilon\in\mathcal{M}_2(\Omega):=\mu-\mu_{|_{A_\varepsilon}}+\lambda \Leb_{A_\varepsilon}\;,\]
i.e.~we replace $\mu$ by $\lambda$ times the Lebesgue measure on the set 
\[A_\varepsilon:=\{x\in\Omega\,:\,\operatorname{dist}(x,\partial\Omega)<\varepsilon\}\;.\]
Let $P:\Omega\to\partial\Omega$ be a map such that $|P(x)-x|=d(x,\partial\Omega)$ for all $x\in A_\varepsilon$. We recall the notation introduced in Subsection \ref{notation}. Then an admissible plan $\bg\in\textsc{Adm}(\mu,\mu_ \varepsilon)$ with $\bg_{\partial\Omega}^{\partial\Omega}=0$ is given by
\begin{equation}\label{eq65}
\begin{aligned}
&\bg_{\Omega}^{ \Omega}=(\operatorname{id},\operatorname{id})_{\sharp}\left(\mu\wedge\mu_ \varepsilon\right)\\
&\bg_{\Omega\setminus A_\varepsilon}^{ A_\varepsilon\cup\partial\Omega}=\bg_{A_\varepsilon\cup\partial\Omega}^{\Omega\setminus A_\varepsilon}=0\\
&\bg_{A_\varepsilon}^{\partial\Omega}=(\operatorname{id},P)_{\sharp}\left(\mu_{|_{A_\varepsilon}}-\lambda\Leb_{A_\varepsilon}\right)_+\\
&\bg_{\partial\Omega}^{A_\varepsilon}=(P,\operatorname{id})_{\sharp}\left(\lambda\Leb_{A_\varepsilon}-\mu_{|_{A_\varepsilon}}\right)_+\,,
\end{aligned}
\end{equation}
where $\mu\wedge\mu_\varepsilon$ denotes the minimum of the measures $\mu$ and $\mu_\varepsilon$, i.e.~their common mass. We seek to show that
\begin{equation}\label{eq66}
\frac{[\mathcal{F}(\mu)-\mathcal{F}(\mu_\varepsilon)]^+}{Wb_2(\mu,\mu_\varepsilon)}=:\frac{J_1(\varepsilon)}{J_2(\varepsilon)}\to+\infty\,,
\end{equation}
for (a sequence of) $\varepsilon\to0$.

We estimate $J_2(\varepsilon)$ using the plan $\bg$ as follows:
\begin{equation}\label{eq:J2}
\begin{aligned}
J_2(\varepsilon)&\le\left(\int_{\overline{\Omega}\times\overline{\Omega}}\,|x-y|^2\,d\bg(x,y)\right)^{1/2}
=\left(\int_{A_\varepsilon}\operatorname{dist}^2(x,\partial\Omega)|\rho(x)-\lambda|\,dx \right)^{1/2}\\
&\le \varepsilon\left(\int_{A_\varepsilon}|\rho(x)-\lambda|\,dx \right)^{1/2}=\varepsilon\,|A_{\varepsilon}|^{1/2}\left(\frac{1}{|A_\varepsilon|}\int_{A_\varepsilon}|\rho(x)-\lambda|\,dx\right)^{1/2}\\
&=:\varepsilon\,|A_{\varepsilon}|^{1/2}M_\varepsilon^{1/2}\,.
\end{aligned}
\end{equation}

Let us now estimate $J_1(\varepsilon)$. Note that there is a constant $c>0$ such that $F(r)\geq c H_\lambda(|r-\lambda|)$ for all $r\geq 0$, where $H_\lambda:[0,\infty)\to[0,\infty)$ is the convex function given by $H_\lambda(s)=F(s+\lambda)$. By the convexity of $F$, noting that $F(\lambda)=0$, and Jensen's inequality we have 
\begin{equation}\label{eq:J1}\begin{aligned}
J_1(\varepsilon)&=\int_{A_\varepsilon}F\left(\rho(x)\right)dx\\
&\geq c  \int_{A_\varepsilon}H_\lambda\left(|\rho(x)-\lambda|\right)\,dx = c |A_\varepsilon| \frac{1}{|A_\varepsilon|}\int_{A_\varepsilon}H_\lambda\left(|\rho(x)-\lambda|\right)\,dx\\
&\geq c |A_\varepsilon| \, H_\lambda(M_\varepsilon)\;.\end{aligned}
\end{equation}
Hence, combining \eqref{eq:J1} and \eqref{eq:J2} we obtain
\[
\frac{J_1(\varepsilon)}{J_2(\varepsilon)}\geq c\frac{|A_\varepsilon|^{1/2}}{\varepsilon}\frac{H_\lambda(M_\varepsilon)}{M_\varepsilon^{1/2}}\;.
\]
Note that $|A_\varepsilon|/\varepsilon\to {\sf H}^{d-1}(\partial\Omega)$ as $\varepsilon\to 0$ and that 
\[\frac{H_\lambda(r)}{\sqrt{r}}>0\quad  \text{ for } r>0\quad \text{and}\quad \lim_{r\to\infty}\frac{H_\lambda(r)}{\sqrt{r}}=+\infty\;.\]
Thus, to obtain \eqref{eq66}, it suffices to show that 
\begin{equation}\label{eq:limsupM}
\limsup_{\varepsilon\to0} M_\varepsilon = \limsup_{\varepsilon\to0} \frac{1}{|A_\varepsilon|}\int_{A_\varepsilon}|\rho(x)-\lambda|dx >0\;.
\end{equation}
This will follow from the assumption that the trace of $G(\rho)$ is not equal to $G(\lambda)$. Note first that this assumption implies that for any $C>\lambda$, we have $G(\rho\wedge C)\in W^{1,2}(\Omega,\Leb_\Omega)$ and $\mathcal T[G(\rho\wedge C)]\neq G(\lambda)$. 
The characterisation of the trace in \cite[Thm. 5.7]{EvGa} yields that
\begin{equation}\label{eq:trG}
\limsup_{\varepsilon\to0} \frac{1}{|A_\varepsilon|}\int_{A_\varepsilon} |G(\rho\wedge C)-G(\lambda)|>0\;.
\end{equation}

Now, if $\lambda >0$, the fact that $G'(\lambda)=\sqrt{\lambda}F''(\lambda)>0$ implies that there is a constant $c>0$ such that  $|r-\lambda|\geq |r\wedge C-\lambda|\geq c|G(r\wedge C)-G(\lambda)|$ for all $r\geq 0$. Hence \eqref{eq:trG} implies the claim \eqref{eq:limsupM}.

If $\lambda=0$, we distinguish two cases, according to whether $0\leq G'(0)=\lim_{r\to0}\sqrt{r}F''(r)$ is finite or not. If $G'(0)<+\infty$, we note that there is a constant $c>0$ such that $r\wedge C\geq c \cdot G(r\wedge C)$ for all $r\geq 0$, and \eqref{eq:limsupM} follows as above. If $G'(0) = +\infty$, we note that we can find a Lipschitz map $\psi:[0,\infty)\to[0,\infty)$ such that $r\wedge C =\psi\big(G(r\wedge C)\big)$ for all $r\geq 0$. We deduce that $\rho\wedge C\in H^1(\Omega)$ and that $\mathcal T[\rho\wedge C]\neq 0$. The characterisation of the trace then yields
\[
\limsup_{\varepsilon\to0} \frac{1}{|A_\varepsilon|}\int_{A_\varepsilon} |\rho\wedge C|>0\;,
\]
which immediately give \eqref{eq:limsupM}. Hence the proof is complete.
\end{proof}

\begin{remark}\label{rem:alternative-proof}
Under the additional assumption that $F$ satisfies the McCann condition, namely $s\mapsto s^dF(s^{-d})$ is convex and non increasing in $(0,+\infty)$, an alternative proof of the lower bound $|\partial \cF|^2(\mu)\geq \cI(\mu)$ is possible following the argument in \cite[Thm.~10.4.6]{AGS} considering directional derivatives of $\cF$. Let us briefly sketch this approach. Given $\mu=\rho\Leb$ with $|\cF(\mu)|<\infty$ and $\xi\in C^\infty_c(\Omega;\R^d)$ and $t$ sufficiently small, consider the perturbation $\mu_t=(\id+t\xi)_\#\mu$. Then one has
\begin{equation}\label{eq:directional}
\lim_{t\searrow 0} \frac{\cF(\mu_t)-\cF(\mu)}{t} = -\int_\Omega L_F(\rho) \nabla\cdot \xi \;d\Leb\;.
\end{equation}
This can be proven verbatim as in \cite[Lem.~10.4.4.]{AGS}.
%Indeed, by the change of variables formula
%\begin{align*}
%\frac{\cF(\mu_t)-\cF(\mu)}{t} = \frac{1}{t}\Big[\int_\Omega F\Big(\frac{\rho(x)}{\det DT_t(x)}\Big)\det DT_t(x)\rho_t(T_tx) dx - \int_\Omega F\big(\rho(x)\big)dx\Big]\\
%= \int_{\Omega} \frac{1}{t}\Big[H\big(\rho(x),\det DT_t(x)\big) - F\big(\rho(x)\big)\Big] dx\;,
%\end{align*}
%where $T_t= \id +t\xi$ and $H(z,s)=sF(z/s)$. Noting that $\partial_s H(z,s)=-L_F(z/s)$ and $d/dt \det DT_t=\nabla\cdot \xi$ and letting $t\searrow 0$, we obtain \eqref{eq:directional}. The passage to the limit under the integral in justified by dominated convergence using the doubling condition on $F$, see the proof of \cite[Lem.~10.4.4.]{AGS}.
Note that $Wb_2(\mu,\mu_t)\leq t\|\xi\|_{L^2(\mu)}$. Now, similar as in the proof of \cite[Thm.~10.4.6 a)]{AGS}, applying \eqref{eq:directional} yields
\[\int_\Omega L_F(\rho) \nabla\cdot \xi \;d\Leb \leq |\partial\cF|(\mu)\|\xi\|_{L^2(\mu)}\;.\]
This allows us to conclude that $L_F(\rho)\in W^{1,1}(\Omega)$ and $\nabla L_F(\rho) = \bw\rho$ with $\|\bw\|_{L^2(\mu)}\leq |\partial F|(\mu)$. From Remark \ref{rem:GvsL} we infer $|\partial \cF(\mu)|^2\geq \cI(\mu)$.
\end{remark}

As an immediate consequence of Proposition \ref{prop4} combined with lower semicontinuity of $\overline{\mathcal{I}}$, see Lemma \ref{lem:lsc-diss}, we get the following
\begin{corollary}\label{cor2}
For all $\mu\in\mathcal{M}_2(\Omega)$ we have $\overline{\mathcal{I}}(\mu)\le|\partial^-\mathcal{F}|^2(\mu)$.
\end{corollary}

The previous bound together with Corollary \ref{cor} yields the following characterisation of solutions to the non-linear diffusion equation \eqref{eq:NLD} with Dirichlet boundary conditions, purely in terms of the function $\cF$ and the metric space $(\mathcal M_2(\Omega),Wb_2)$.

\begin{theorem}\label{theo3}
The relaxed slope $|\partial^-\cF|$ is a strong upper gradient of $\cF$ on $(\mathcal M_2(\Omega),Wb_2)$. A curve $(\mu_t)_{t\in[0,T]}$ is a curve of maximal slope w.r.t. this strong upper gradient if and only if $\mu_t=\rho_t\Leb_\Omega$ with $t\mapsto G(\rho_t)-G(\lambda)$ belonging to $L^2([0,T];W^{1,2}_0(\Omega))$ such that $(\rho_t)$ is a weak solution to \eqref{eq:NLD}.\end{theorem}

\begin{proof}
Since $\sqrt{\overline{\mathcal{I}}}$ is a strong upper gradient, so it is $|\partial^-\cF|$ by the bound in Corollary \ref{cor2}. For the same reason, any curve of maximal slope w.r.t. $|\partial^-\cF|$ is also a curve of maximal slope w.r.t. $\sqrt{\overline{\mathcal{I}}}$. Hence the conclusion follows from Corollary \ref{cor} and in turn from Theorem \ref{teo2}.
\end{proof}

\subsection{Convergence of the minimising movement scheme}

In this section we consider the minimising movement scheme for the metric gradient flow, a time-discrete variational approximation scheme for curves of maximal slope. In the setting of gradient flows in the space of probability measures, this scheme is also known as the Jordan–Kinderlehrer–Otto scheme \cite{JKO}. Leveraging abstract results on minimising movements together with the characterisation of curves of maximal slope in Theorem \ref{theo3}, we obtain that the scheme converges in the present setting to solutions of the non-linear diffusion with Dirichlet boundary conditions. This partially extends the results obtained in \cite{FG10, KKS22} to more general non-linear equations.
\medskip

Let us fix an initial measure $\mu_0\in \mathcal{M}_2(\Omega)$ such that $\mathcal{F}(\mu_0)<\infty$ (the assumption of finiteness of the entropy at the initial point ensures that the measure is absolutely continuous given the lack of convexity of $\mathcal{F}$). Given a time step $\tau>0$, we consider the recursive sequence $(\mu_n^\tau)_n\in\mathcal{M}_2(\Omega)$ defined as
\begin{equation}\label{eq:scheme}
\mu_0^\tau:=\mu_0,\quad\text{and}\quad\mu_{n+1}^\tau:=\argmin_\nu \left\{\mathcal{F}(\nu)+\frac{Wb_2(\nu,\mu_n^\tau)}{2\tau}\right\}\,.
\end{equation}
Then we build a discrete gradient flow trajectory as a piecewise constant interpolation $(\bar\mu_t^\tau)_{t\ge0}\in\mathcal{M}_2(\Omega)$ given by 
$$
\bar\mu_0^\tau:=\mu_0,\quad\text{and}\quad\bar\mu^\tau(t):=\bar\mu_n^\tau\quad\text{for}\,\,\,t\in[n\tau,(n+1)\tau)\,.
$$
Then we have the following result
\begin{theorem}\label{JKO}
For any $\tau>0$ and $\mu_0\in\mathcal{M}_2(\Omega)$ with $\mathcal{F}(\mu_0)<\infty$ the variational scheme \eqref{eq:scheme} admits a solution $(\mu_n^\tau)_n$. As $\tau\to0$, for any family of discrete solutions there exists a sequence $\tau_k\to0$ and a $2$-absolutely continuous curve $(\mu_t)_{t\ge0}$ such that
$\bar\mu^{\tau_k}_t$ converges to some limit measure $\mu_t$ in $(\mathcal{M}_2(\Omega), Wb_2)$ as $k\to\infty$ for all $t\in[0,+\infty)$. Moreover, any such limit curve is of the form $\mu_t=\rho_t\Leb_\Omega$ with $t\mapsto G(\rho_t)-G(\lambda)$ belonging to $L^2([0,T];W^{1,2}_0(\Omega))$ such $(\rho_t)$ is a weak solution to the diffusion equation \eqref{eq:NLD}.
\end{theorem}

\begin{proof}
The result follows from general results for metric gradient flows, see \cite[Section 2.3]{AGS} for a deep analysis of the so called minimizing movement scheme. We consider the metric space $(\mathcal{M}_2(\Omega), Wb_2)$ endowed with the vague topology. It follows that \cite[Assumptions 2.1 (a,b,c)]{AGS} are satisfied. Existence of a solution to the variational scheme \eqref{eq:scheme} and of a subsequential limit curve $(\mu_t)_t$ now follows from \cite[Corollary 2.2.2, Proposition 2.2.3]{AGS}. Moreover, \cite[Theorem 2.3.1]{AGS} gives that the limit curve is a curve of maximal slope for the strong upper gradient $|\partial^-\mathcal F|$, see \eqref{eq:relslope}. Thus, by Corollary \ref{cor2} it is also a curve of maximal slope with respect to the strong upper gradient $\sqrt{\overline{\mathcal{I}}}$. Theorem \ref{teo2} yields the identification with a weak solution to \eqref{eq:NLD}.
\end{proof}

\appendix

\section{A boundary estimate}

Here we provide a technical Lemma to estimate the $L^2$ norm of Sobolev functions with vanishing trace in tubular neighborhoods of the boundary. Similar results are classically proven in balls around boundary points \cite[Thm.~5.7]{EvGa}, see also a version with very narrow neighborhood in (the proof of ) \cite[Lemma A.3]{CMS}.\smallskip

Let $\Omega\subset \R^n$ be a bounded open domain with $C^1$ boundary and set $B_\kappa:=\{x\in \Omega: {\rm dist}(x,\partial\Omega)<\kappa\}$ for $\kappa>0$.

\begin{lemma}\label{lem:bdry-est}
There exists a constant $C_\Omega>0$ such that for any $f\in W^{1,2}_0(\Omega)$ and $\kappa$ sufficiently small, we have
\begin{equation*}
\int_{B_\kappa} |f|^2 d\Leb \leq C_\Omega \kappa^2 \int_{B_\kappa} |\nabla f|^2 d\Leb\;. 
\end{equation*}
\end{lemma}

\begin{proof}
By a standard argument using a partition of unity and straightening the boundary, we can assume without restriction that $f$ is supported in a ball $U=B_r(\bar x)$ with $\bar x\in \partial \Omega$ and $\Omega$ is the halfspace $\{x\in \R^n\colon x_1>0\}$. By density in $W^{1,2}_0(\Omega)$, it is further sufficient to consider $f$ that belongs to $C^\infty_c(\Omega)$. Writing $x=(r,x')\in \R^+\times\R^{n-1}$, and using that $f(0,x')=0$ for all $x'$, we now estimate
\begin{align*}
\int_{B_\kappa} |f(x)|^2 dx 
&= 
\int_{\R^{n-1}} \int_0^\kappa |f(r,x')|^2\,dr d x'
=
\int_{\R^{n-1}} \int_0^\kappa\left|\int_0^{r}\partial_1 f(s,x')\,ds\right|^2\,drdx'\\
&\leq
\int_{\R^{n-1}} \int_0^\kappa r\int_0^{r}\left|\partial_1 f(s,x')\right|^2\,ds drdx'
\leq
\frac{\kappa^2}{2}\int_{\R^{n-1}} \int_0^\kappa |\partial_1 f(s,x')|^2\,dsdx'\\
&\leq
\frac{\kappa^2}{2}\int_{B_\kappa}|\nabla f|^2d\Leb\;.
\end{align*}
\end{proof}

\section{Existence of weak solutions}

Here we provide a simple existence result for weak solutions to the non-linear diffusion equation \eqref{eq:NLD} with Dirichlet boundary conditions.

\begin{lemma}\label{lem:weaksol-energy}
For any $\mu_0=\rho_0\Leb_\Omega$ with $\cF(\mu_0)<\infty$ there exists a weak solution $(\rho_t)$ to the non-linear diffusion equation \eqref{eq:NLD} with initial datum $\rho_0$ such that $t\mapsto G(\rho_t)-G(\lambda)$ belongs to $L^2\big([0,T],W^{1,2}_0(\Omega)\big)$ for any $T>0$.
\end{lemma}

\begin{proof}
We slightly adapt the classical argument for the existence of energy solutions to the generalised porous medium equation, see \cite[Sec. 5.4]{Vaz07}.
\medskip

\emph{Step1:} Recall that $\cF(\mu_0)<\infty$ implies that $\rho_0\in L^1(\Omega)$.
Approximate the domain $\Omega$ by an increasing sequence of domains $\Omega_n\subset \Omega$ with smooth boundary and let $B_n$ be the $1/n$-neighbourhood of $\Omega\setminus \Omega_n$ in $\Omega$. Set $\lambda_n:=\min\big\{\max\big\{1/n, \lambda\}, n\big\}$ and consider a sequence of approximate initial data $\rho^n_0$ such that 
\begin{equation}\label{eq:boundsrho}
1/n \leq \rho_0^n\leq n\;, \quad \rho_0^n = \lambda_n \text{ on } B_n\;, \quad \rho_0^n\to \rho_0 \text{ in } L^1(\Omega)\;, \quad \cF(\mu_0^n)\to \cF(\mu_0)
\end{equation}
as $n\to\infty$ for $\mu_0^n=\rho_0^n\Leb_\Omega$.
This can be achieved e.g. by setting
\[
\rho_0^n = \min\big\{\max\big\{1/n, \rho_0\}, n\big\}\cdot \boldsymbol{1}_{\Omega\setminus B_n} + \lambda_n\cdot \boldsymbol{1}_{B_n}\;.
\]
Now, there exists a unique classical solution $\rho^n\in C^{1,2}\big((0,\infty)\times \overline{\Omega}\big)$ to the approximate problem
$$
\begin{cases}
\partial_t \rho^n = \Delta L_F(\rho^n)\quad &\text{ in } (0,\infty)\times \Omega_n\;,\\
\rho^n(0,\cdot ) = \rho^n\quad\quad\quad\quad &\text{ in } \overline\Omega_n\;,\\
\rho^n = \lambda_n\qquad\quad\quad &\text{ on } [0,\infty)\times\partial \Omega_n\;,
\end{cases}
$$
which satisfies $1/n\leq \rho^n_t\leq n$ for all $t\geq 0$ and which attains the initial datum in the sense that $\rho^n_t\to\rho^n_0$ in $L^1(\Omega_n)$ as $t\to0$. This solution is obtained by classical existence results and maximum principles for non-degenerate quasilinear equations, noting that the potential degeneracy of  $L_F'$ at $0$ or $\infty$ is avoided due to the upper and lower bounds. We extend $\rho^n_t$ to $\Omega$ by the constant value $\lambda_n$. 
Let us denote by $F_n$ the modification of $F$ as in Remark \ref{rem:mod-F}, i.e. $F_n:[0,+\infty)\to[0,+\infty)$
$$
F_n(r)=F(r)-F(\lambda_n)-F'(\lambda_n)(r-\lambda_n)\,.
$$ 
It takes its unique minimum at $\lambda_n$ and recall that $F_n''=F''$. Let $\cF_n$ be the corresponding modification of $\cF$ replacing $F$ in Definition \ref{def:energy} by $F_n$. Note that 
\begin{equation}\label{eq:FnvsF}
\big|\cF_n(\mu)-\cF(\mu)\big| \leq F(\lambda_n) |\Omega|+ |F'(\lambda_n)| \cdot\|\rho-\lambda_n\|_{L^1(\Omega)}\qquad \forall \mu=\rho\Leb_\Omega\;.
\end{equation}

Setting $\mu^n_t=\rho^n_t\Leb_\Omega\in \mathcal{M}_2(\Omega)$, due to the regularity and bounds for $\rho^n$ we obtain by direct calculation for any $t>0$:
\begin{align*}
\frac{d}{dt}\mathcal{F}_n(\mu^n_t) 
&= 
\int_\Omega F_n'(\rho^n_t)\partial_t\rho^n_t\; d\Leb
 = 
\int_{\Omega_n} F_n'(\rho_t) \Delta L_F(\rho^n_t)\; d\Leb\\
&= 
-\int_{\Omega_n} |F''(\rho^n_t)|^2 \rho^n_t |\nabla\rho^n_t|^2\; d\Leb
=
-\int_\Omega |\nabla G(\rho^n_t)|^2\; d\Leb
=
-\mathcal I(\mu^n_t)\;.
\end{align*}
Here, the integration domain $\Omega$ can be exchanged with $\Omega_n$, since $\rho^n_t$ is constant equal to $\lambda_n$ on a neighbourhood of $\Omega\setminus \Omega_n$. No boundary term occurs in the integration by parts, since $F_n'(\rho^n_t)=F_n'(\lambda_n)=0$ on $\partial \Omega_n$.
Note that $\cF_n(\mu^n_t)\to \cF_n(\mu_0^n)$ as $t\to0$ by dominated convergence due to the convergence of $\rho^n_t$ in $L^1$ and the uniform upper and lower bounds on $\rho^n_t$. Thus, integrating the last display on $[0,T]$ we obtain
\begin{equation}\label{eq:energy-id-n}
\cF_n(\mu^n_T) - \cF_n(\mu^n_0) +\int_0^T\cI(\mu^n_t)dt = 0\;.
\end{equation}
From \eqref{eq:boundsrho} and \eqref{eq:FnvsF} we infer that $\cF_n(\mu^n_0)\to \cF(\mu_0)<\infty$. Since $\cF_n\geq 0$, we deduce that $\nabla G(\rho^n)$ is bounded uniformly in $n$ in $L^2\big((0,T)\times\Omega\big)$ and by Poincar\'e inequality so is $G(\rho^n)-G(\lambda_n)$.\medskip

\emph{Step 2:} The solutions obtained in the previous step are contractive in $L^1(\Omega_n)$, see \cite[Sec.~3.2.3]{Vaz07}. Hence we infer that for $n\leq m$ and all $t\geq 0$ we have
$\|\rho^n_t-\rho^m_t\|_{L^1(\Omega)} \leq \|\rho^n_0-\rho^m_0\|_{L^1(\Omega)}+\Leb(B_n)/n$. Thus, for any $T>0$ we conclude that as $n\to\infty$  the curve $\rho^n$ converges to a limit $\rho$ in $C\left([0,T];L^1(\Omega)\right)$. Moreover, from the uniform boundedness in $L^2$ we have that up to subsequences  $G(\rho^n)-G(\lambda_n)$ and $\nabla G(\rho^n)$ have strong respectively weak limits in $L^2\left((0,T)\times\Omega\right)$. This implies that $G(\rho)-G(\lambda)$ belongs to $L^2\big((0,T);W^{1,2}_0(\Omega)\big)$ that it is the limit of the whole sequence $G(\rho^n)-G(\lambda_n)$ in this space. Continuity of $t\mapsto \rho_t$ in $L^1(\Omega)$ trivially implies vague continuity of the curve. It remains to verify that $\rho$ is a weak solution of \eqref{eq:NLD}. For the classical solution $\rho^n$ we have for any $\varphi\in C^\infty_c(\Omega_n)$ and $0\leq s<t$:
\[
\int_\Omega\varphi (\rho^n_t-\rho^n_s)\,d\Leb = \int_s^t \int_\Omega \Delta \varphi\, L_{F}(\rho^n_r) \,d\Leb\,dr = -\int_s^t\int_\Omega \langle\nabla\varphi,\nabla L_F(\rho^n_r)\rangle\, d\Leb\, dr\;.
\]
Using that $\nabla L_F(\rho^n)=\sqrt{\rho^n}\cdot\nabla G(\rho^n)$, the convergence of $\rho^n$ to $\rho$ in $C\big([0,T];L^1(\Omega)\big)$, and the weak convergence of $\nabla G(\rho^n)$ to $\nabla G(\rho)$ in $L^2\big((0,T)\times\Omega\big)$, we can pass to the limit in the previous equation and obtain \eqref{eq:weak-sol} after a further integration by parts.
\end{proof}

\begin{remark}\label{rem:weak-sol}
The variational characterisation of weak solutions in Theorem \ref{teo2}, shows that \eqref{eq:energy-id-n} passes to the limit and the solution constructed above satisfies the energy-dissipation identity
\[
\cF(\mu_T) - \cF(\mu_0)= -\int_0^T\cI(\mu_t)dt \leq 0\;.
\]
\end{remark}

\end{document}